\documentclass [12pt]{article}
\usepackage{amssymb}
\def\lanbox{\hbox{$\, \vrule height 0.25cm width 0.25cm depth 0.01cm \,$}}

\usepackage{color}

\begin{document}

\centerline{\Large\bf Solvability conditions for some non-Fredholm}

\centerline{\Large\bf operators with shifted arguments}

\bigskip

\centerline{Vitali Vougalter$^{1 \ *}$, Vitaly Volpert$^{2, 3}$ }

\bigskip

\centerline{$^{1 \ *}$ Department of Mathematics, University
of Toronto, Toronto, Ontario, M5S 2E4, Canada}

\centerline{e-mail: vitali@math.toronto.edu}

\bigskip

\centerline{$^2$ Institute Camille Jordan, UMR 5208 CNRS,
University Lyon 1, Villeurbanne, 69622, France}

\centerline{$^3$ Peoples' Friendship University of Russia, 6 Miklukho-Maklaya
St,}

\centerline{Moscow, 117198, Russia}

\centerline{e-mail: volpert@math.univ-lyon1.fr}

\bigskip
\bigskip
\bigskip

\noindent {\bf Abstract.}
In the first part of the article we establish the existence in the sense of sequences of solutions in $H^{2}({\mathbb R})$
for some nonhomogeneous linear differential equation in which one of the terms has the argument translated by a constant.
It is shown that under the reasonable technical conditions the convergence in
$L^{2}({\mathbb R})$ of the source terms implies the existence and the convergence
in $H^{2}({\mathbb R})$ of the solutions. The second part of the work deals with the solvability in the sense of sequences
in $H^{2}({\mathbb R})$ of the integro-differential equation in which one of the terms has the argument shifted by a constant.
It is demonstrated that under the appropriate auxiliary assumptions the convergence in
$L^{1}({\mathbb R})$ of the integral kernels yields the existence and the convergence
in $H^{2}({\mathbb R})$ of the solutions.
Both equations considered involve the second order differential operator with or without the Fredholm
property depending on the value of the constant by which the argument gets translated.

\bigskip
\bigskip



\bigskip
\bigskip
\bigskip
\bigskip

\section{Introduction}


 Let us recall that a linear operator $L$ acting from a Banach
 space $E$ into another Banach space $F$ possesses the Fredholm
 property if its image is closed, the dimension of its kernel and
 the codimension of its image are finite. As a consequence, the
 problem $Lu=f$ is solvable if and only if $\phi_i(f)=0$ for a
 finite number of functionals $\phi_i$ from the dual space $F^*$.
 These properties of the Fredholm operators are broadly used in different
 approaches of the linear and nonlinear analysis.

 Elliptic problems in bounded domains with a sufficiently smooth
 boundary satisfy the Fredholm property if the ellipticity
 condition, proper ellipticity and Lopatinskii conditions are
 fulfilled (see e.g. \cite{Ag}, \cite{E09}, \cite{LM}, \cite{Volevich}).
 This is the main result of the theory of linear elliptic equations.
 When domains are unbounded, these conditions may
 be insufficient and the Fredholm property may not be satisfied.
 For example, the Laplace operator, $Lu = \Delta u$, in $\mathbb R^d$
 does not satisfy the Fredholm property when considered in
 H\"older spaces, $L : C^{2+\alpha}(\mathbb R^d) \to C^{\alpha}(\mathbb
 R^d)$, or in Sobolev spaces,  $L : H^2(\mathbb R^d) \to L^2(\mathbb
 R^d)$.

 Linear elliptic problems in unbounded domains satisfy the
 Fredholm property if and only if, in addition to the conditions
 given above, the limiting operators are invertible (see \cite{V11}). In some
 trivial situations, the limiting operators can be constructed explicitly . For
 example, if

 $$ L u = a(x) u'' + b(x) u' + c(x) u , \;\;\; x \in \mathbb R,
 $$
 with the coefficients of this operator having limits at infinities,

 $$ a_\pm =\lim_{x \to \pm \infty} a(x) , \;\;\;
 b_\pm =\lim_{x \to \pm \infty} b(x) , \;\;\;
 c_\pm =\lim_{x \to \pm \infty} c(x) , $$
 the limiting operators are given by:

 $$ L_{\pm}u = a_\pm u'' + b_\pm u' + c_\pm u . $$
 Because the coefficients here are constants, the essential spectrum of
 the operator, that is the set of complex numbers $\lambda$ for
 which the operator $L-\lambda$ does not satisfy the Fredholm
 property, can be explicitly computed by virtue of the Fourier
 transform (\ref{f}):

 $$ \lambda_{\pm}(\xi) = -a_\pm \xi^2 + b_\pm i\xi + c_\pm , \;\;\;
 \xi \in \mathbb R . $$
 The invertibility of the limiting operators is equivalent to the
 condition that the origin does not belong to the essential spectrum.

 In the case of the general elliptic problems, the same assertions hold true.
 The Fredholm property is satisfied if the essential spectrum does not contain
 the origin or if the limiting operators are invertible. However,
 such conditions may not be written explicitly.

 When the operators are  non-Fredholm, the usual solvability
 conditions may not be applicable and the solvability relations
 are, in general, unknown. There are certain classes of operators
 for which the solvability conditions are obtained. Let us illustrate
 them with the following example. We consider the equation

 \begin{equation}
 \label{int1}
 L_{0}u \equiv -\Delta u -a u = f
 \end{equation}
 in $\mathbb R^d$, where $a>0$ is a constant.
 The self-adjoint operator $L_{0}$ here coincides with its limiting operators.
 The homogeneous problem admits a nonzero bounded solution. The Sobolev space
 used throughout the article is equipped with the norm
 \begin{equation}
 \label{sob}
 \|u\|_{H^{2}(\mathbb R^{d})}^{2}:=\|u\|_{L^{2}(\mathbb R^{d})}^{2}+
 \|\Delta u\|_{L^{2}(\mathbb R^{d})}^{2}, \quad d\in {\mathbb N}.
 \end{equation}
 Note that for the
 Fredholm property we consider bounded but not $H^2(\mathbb R^d)$ solutions
 of the corresponding homogeneous adjoint problem. Because equation (\ref{int1})
 has constant coefficients, we can apply the standard Fourier transform
 to solve it explicitly. If $f \in L^2(\mathbb R^d)$ and
 $xf \in L^1(\mathbb R^d)$, then (\ref{int1}) admits
 a unique solution in $H^2(\mathbb R^d)$ if and
 only if
 \begin{equation}
 \label{oc0}
 \Bigg(f(x),\frac{e^{ipx}}{(2\pi)^{\frac{d}{2}}}\Bigg)_{L^2(\mathbb R^d)}=0, \quad
 p\in S_{\sqrt{a}}^{d} \quad a.e.
 \end{equation}
 (see Lemmas 5 and 6 of \cite{VV103}).
We denote
\begin{equation}
\label{ip}
(f_{1}(x),f_{2}(x))_{L^{2}({\mathbb R}^{d})}:=\int_{{\mathbb R}^{d}}f_{1}(x)
\bar{f_{2}}(x)dx,
\end{equation}
with a slight abuse of notations when these functions are not square
integrable. Indeed, if $f_{1}(x)\in L^{1}({\mathbb R}^{d})$ and $f_{2}(x)$ is
bounded, then the integral in the right side of definition (\ref{ip}) makes
sense. The orthogonality relations (\ref{oc0}) are with respect to the
standard Fourier harmonics solving the homogeneous adjoint problem for
(\ref{int1}), belonging to $L^{\infty}(\mathbb R^d)$ but they are not square
 integrable.

 $S_{\sqrt{a}}^{d}$ in (\ref{oc0}) stands for
 the sphere in $\mathbb R^d$ of radius $\sqrt{a}$ centered at the origin.
 Hence, though our operator does not satisfy the Fredholm property, the
 solvability conditions are formulated in a similar way. However,
 this similarity is only formal because the range of the operator
 is not closed.

 In the case of the operator involving a scalar potential,
 $$ L u \equiv -\Delta u -a(x) u = f ,
 $$
 the standard Fourier transform is not of any help. Nevertheless, the solvability
 relations in ${\mathbb R}^{3}$ can be obtained by a rather sophisticated
 application of the theory of self-adjoint
 operators (see \cite{VV08}). As before, the solvability conditions are
 formulated in terms of the orthogonality to the solutions of the homogeneous adjoint
 problem. There are several other examples of linear elliptic
 operators without the Fredholm property for which the solvability
 relations can be derived (see ~\cite{EV21}, ~\cite{V11}, ~\cite{VKMP02}, ~\cite{VV1031},
 ~\cite{VV08}, ~\cite{VV10}, ~\cite{VV103}).

 Solvability conditions play a significant role in the analysis of
 nonlinear elliptic problems. When the operators do not satisfy the Fredholm
 property, in spite of a certain progress in understanding the linear
 equations, there exist only few examples where nonlinear non-Fredholm
 operators are analyzed (see ~\cite{DMV05}, ~\cite{DMV08}, ~\cite{DMV08b},
 ~\cite{EV211}, ~\cite{EV25}, ~\cite{VV1021}, ~\cite{VV103}, ~\cite{VV25}).
 The operators without the Fredholm
 property arise also when studying the so-called embedded solitons
 (see e.g. ~\cite{PY02}). Wave systems with an infinite number of localized traveling
waves were covered in ~\cite{AMP14}.
Standing lattice solitons in the discrete NLS equation with saturation
were discussed in ~\cite{AKLP19}.
Fredholm structures,
 topological invariants and applications were considered in ~\cite{E09}.
Book ~\cite{E18} deals with the symmetrization and stabilization of solutions of nonlinear elliptic equations.
 Exponential decay of solutions, Fredholm and properness properties of
 quasilinear elliptic systems of
 second order were studied in ~\cite{GS05} and ~\cite{GS10}. The existence of solutions for the
boundary value problems for differential-difference equations with incommensurable shifts was established in ~\cite{S92}.
The existence of classical solutions of hyperbolic equations with translation operators in lower-order derivatives was covered
in ~\cite{Z25}. A classical solution to a hyperbolic differential-difference equation with a translation by an arbitrary vector was
constructed in ~\cite{ZM23}. Elliptic equations with translations of general form in a half-space were discussed in ~\cite{M22}.
In work ~\cite{MZ23} the authors construct the classical solutions of hyperbolic differential-difference equations with differently directed translations. Elliptic problems with shift arise in the investigation of travelling waves in various models in mathematical biology with time delay (see \cite{TV}).

In the first part of the present article we consider the linear problem
\begin{equation}
\label{eqsh1}
-\frac{d^{2}u(x)}{dx^{2}}-au(x-h)=f(x), \quad x\in {\mathbb R}
\end{equation}
with constants $a>0, \ h\in {\mathbb R}, \ h\neq 0$ and the square
integrable right side. Note that equation (\ref{eqsh1}) is the one dimensional
analog of (\ref{int1}) obtained by shifting the argument by $h$ in the second
term of the left side. The non self-adjoint operator
\begin{equation}
\label{Lh}
L_{h}: H^{2}({\mathbb R})\to L^{2}({\mathbb R})
\end{equation}
involved in the left side of
(\ref{eqsh1}) has the essential spectrum $\lambda_{h}(p)$ which can be easily
found using the standard Fourier transform (\ref{f}). Hence,
\begin{equation}
\label{lhp}
\lambda_{h}(p)=p^{2}-a \cos(ph)+ia \sin(ph),
\end{equation}
so that
\begin{equation}
\label{lhp2}
|\lambda_{h}(p)|^{2}=(p^{2}-a)^{2}+2ap^{2}(1- \cos(ph))
\end{equation}
with $p\in {\mathbb R}$ and the constants $a>0, \ h\in {\mathbb R}, \ h\neq 0$.
Our argument shift by $h$ is a unitary operator on $L^{2}({\mathbb R})$.
For the Fourier transform and the norm we have
\begin{equation}
\label{uxhpn}
{\cal F}[u(x-h)](p)=\widehat{u}(p)e^{-iph}, \quad
\|u(x-h)\|_{L^{2}({\mathbb R})}=\|u(x)\|_{L^{2}({\mathbb R})},
\end{equation}
where $u(x)\in L^{2}({\mathbb R})$.

Evidently, if
$\displaystyle{h\neq \frac{2\pi n}{\sqrt{a}}, \ n\in {\mathbb Z}}$, as
distinct from the situation when the shift constant $h$ vanishes, our operator
$L_{h}$ satisfies the Fredholm property since its essential spectrum does not
contain the origin. In this case we have
\begin{equation}
\label{alfa}
|\lambda_{h}(p)|^{2}=(p^{2}-a)^{2}+2ap^{2}(1- \cos(ph))\geq \alpha>0, \quad
p\in {\mathbb R}
\end{equation}
for a certain positive constant $\alpha$.

On the other hand, if
$\displaystyle{h=\frac{2\pi n}{\sqrt{a}}, \ n\in {\mathbb Z}, \ n\neq 0}$, the
operator $L_{h}$ fails to satisfy the Fredholm property because the origin
belongs to its essential spectrum, which can be easily seen from formula
(\ref{lhp2}). The solvability of equation (\ref{eqsh1}) is covered in Theorem 1
below.

The second part of the present work deals with the nonlinear problem, for
which the Fredholm property may or may not be satisfied:
\begin{equation}
\label{id1}
\frac{d^{2}u(x)}{dx^{2}}+au(x-h)+\int_{-\infty}^{\infty}G(x-y)F(u(y),y)dy=0, \quad x\in {\mathbb R}
\end{equation}
with constants $a>0, \ h\in {\mathbb R}, \ h\neq 0$. The argument in the
second term in the left side of (\ref{id1}) is shifted by $h$ as in
equation (\ref{eqsh1}). The solvability of the equation analogous to
(\ref{id1}) when $h$ vanishes was considered in ~\cite{VV1021} via the
fixed point technique.
The integro-differential equations appear in cell population dynamics in the
situations when there is a nonlocal consumption of resources and intraspecific
competition. The space variable $x$ in such problems
corresponds to the cell genotype and not to the usual physical space.
The function $u(x)$ stands for the cell density as a
function of the genotype. The evolution of cell density in the time dependent
problems is due to the cell proliferation and mutations.
The diffusion term corresponds to the change of genotype via small random
mutations, and the integral term describes large mutations. In such context,
$F(u,x)$ is the rate of cell birth which depends on $u$ and $x$ (density
dependent proliferation), and the function $G(x-y)$ shows the proportion of
newly born cells changing their genotype from $y$ to $x$. Let us assume that
it depends on the distance between the genotypes.
In the population dynamics the
integro-differential equations describing models with intra-specific
competition and nonlocal consumption of resources were studies actively
in recent years (see e.g.
~\cite{ABVV10}, ~\cite{BNPR09}, ~\cite{GVA06}).


\setcounter{equation}{0}

\section{Formulation of the results}

Our first main statement is as follows.

\bigskip

\noindent
{\bf Theorem 1.} {\it Let $f(x)\in L^{2}({\mathbb R})$.

\noindent
a) If $\displaystyle{h\neq \frac{2\pi n}{\sqrt{a}}, \ n\in {\mathbb Z}}$ then
equation (\ref{eqsh1}) admits a unique solution $u(x)\in H^{2}({\mathbb R})$.

\noindent
b) When $\displaystyle{h=\frac{2\pi n}{\sqrt{a}}, \ n\in {\mathbb Z}, \ n\neq 0}$,
let in addition $xf(x)\in L^{1}({\mathbb R})$. Then problem
(\ref{eqsh1}) possesses a unique solution $u(x)\in H^{2}({\mathbb R})$ if
and only if the orthogonality conditions
\begin{equation}
\label{ocf0}
\Bigg(f(x), \frac{e^{\pm i\sqrt{a}x}}{\sqrt{2\pi}}\Bigg)_{L^{2}({\mathbb R})}=0
\end{equation}
hold.}

Note that in the situation b) of the theorem above the solvability conditions
of equation (\ref{eqsh1}) are analogous to the ones presented in the case a)
of Lemma 5 of ~\cite{VV103} when the shift parameter $h$ was trivial.

We introduce the sequence of linear approximate equations related to problem
(\ref{eqsh1}) as
\begin{equation}
\label{eqsh1m}
-\frac{d^{2}u_{m}(x)}{dx^{2}}-au_{m}(x-h)=f_{m}(x), \quad x\in {\mathbb R}, \quad
m\in {\mathbb N}
\end{equation}
with the constants $a>0, \ h\in {\mathbb R}, \ h\neq 0$.  The sequence of functions
$\{f_{m}(x)\}_{m=1}^{\infty}$ converges to $f(x)$ as $m\to \infty$ in the appropriate
function spaces given below. We demonstrate that, under the reasonable assumptions,
each equation (\ref{eqsh1m}) possesses a unique solution $u_{m}(x)\in H^{2}({\mathbb R})$,
limiting problem  (\ref{eqsh1})  has a unique solution $u(x)\in H^{2}({\mathbb R})$ and
$u_{m}(x)\to u(x)$ in $H^{2}({\mathbb R})$ as $m\to \infty$. This is the so-called  {\it existence of solutions in the sense
of sequences} (see ~\cite{V11}, ~\cite{VV1031}).
Our second main result
is as follows.

\bigskip

\noindent
{\bf Theorem 2.} {\it Let $m\in {\mathbb N}, \ f_{m}(x)\in L^{2}({\mathbb R})$,
such that $f_{m}(x)\to f(x)$ in $L^{2}({\mathbb R})$ as $m\to \infty$.

\noindent
a) If $\displaystyle{h\neq \frac{2\pi n}{\sqrt{a}}, \ n\in {\mathbb Z}}$ then
each equation (\ref{eqsh1m}) has a unique solution
$u_{m}(x)\in H^{2}({\mathbb R})$ and limiting problem (\ref{eqsh1}) admits a
unique solution $u(x)\in H^{2}({\mathbb R})$.

\noindent
b) When $\displaystyle{h=\frac{2\pi n}{\sqrt{a}}, \ n\in {\mathbb Z}, \ n\neq 0}$,
let in addition $xf_{m}(x)\in L^{1}({\mathbb R})$, such that
$xf_{m}(x)\to xf(x)$ in $L^{1}({\mathbb R})$ as $m\to \infty$.
Moreover, we assume that the orthogonality relations
\begin{equation}
\label{ocfm}
\Bigg(f_{m}(x), \frac{e^{\pm i\sqrt{a}x}}{\sqrt{2\pi}}\Bigg)_{L^{2}({\mathbb R})}=0
\end{equation}
are valid for all $m\in {\mathbb N}$. Then each problem (\ref{eqsh1m}) admits a
unique solution $u_{m}(x)\in H^{2}({\mathbb R})$ and limiting equation
(\ref{eqsh1}) possesses a unique solution $u(x)\in H^{2}({\mathbb R})$.

\noindent
In both cases a) and b), we have $u_{m}(x)\to u(x)$ in $H^{2}({\mathbb R})$
as $m\to \infty$.}

The nonlinear part of problem (\ref{id1}) will satisfy the
regularity conditions analogous to the ones of ~\cite{VV1021}.

\bigskip

\noindent
{\bf Assumption 3.} {\it Function $F(u,x): {\mathbb R}\times {\mathbb R}
\to {\mathbb R}$ is satisfying the Caratheodory condition (see ~\cite{K64}), such that
\begin{equation}
\label{ub1}
|F(u,x)|\leq k|u|+h(x) \quad for \quad u\in {\mathbb R}, \ x\in {\mathbb R},
\end{equation}
with a constant $k>0$ and $h(x):{\mathbb R}\to {\mathbb R}^{+}, \quad
h(x)\in L^{2}({\mathbb R})$. Furthermore, it is a Lipschitz continuous function,
such that
\begin{equation}
\label{lk1}
|F(u_{1},x)-F(u_{2},x)|\leq l |u_{1}-u_{2}| \quad for \quad any \quad
 u_{1,2}\in{\mathbb R}, \quad x\in {\mathbb R}
\end{equation}
with a constant $l>0$.}

\bigskip

 Article  ~\cite{BO86} is devoted to the solvability of a local elliptic problem in a bounded domain
 in ${\mathbb R}^{N}$. The nonlinear function involved there was allowed to have a sublinear growth.
In order to establish the solvabity of equation (\ref{id1}), we will use the auxiliary problem
\begin{equation}
\label{auxnl}
-\frac{d^{2}u(x)}{dx^{2}}-au(x-h)=\int_{-\infty}^{\infty}G(x-y)F(v(y), y)dy, \quad x\in {\mathbb R},
\end{equation}
where $a>0, \ h\in {\mathbb R}, \ h\neq 0$ are the constants.
The main issue for equation (\ref{auxnl}) is that operator (\ref{Lh}) contained in its left side may not
satisfy the Fredholm property depending on the value of the constant $h$ as discussed above, which is the
obstacle when solving such problem. Similar cases
appearing in linear and nonlinear problems, both self- adjoint and non
self-adjoint containg non-Fredholm second or fourth order differential
operators or even systems of equations with non-Fredholm operators
have been treated extensively in recent years (see ~\cite{EV21}, ~\cite{EV211}, ~\cite{EV25}, ~\cite{V11}, ~\cite{VKMP02}, ~\cite{VV1031},
~\cite{VV08}, ~\cite{VV10}, ~\cite{VV1021}, ~\cite{VV103}, ~\cite{VV25}). Our next main proposition is as follows.

\bigskip

\noindent
{\bf Theorem 4.}  {\it Let  $G(x): {\mathbb R}\to {\mathbb R}, \ G(x)\in L^{1}({\mathbb R})$, the constants
$a>0, \ h\in {\mathbb R}, \ h\neq 0$
and Assumption 3  is valid.

\noindent
a) If $\displaystyle{h\neq \frac{2\pi n}{\sqrt{a}}, \ n\in {\mathbb Z}}$, we assume that $2\sqrt{\pi}N_{a, \ h}l<1$ with
$N_{a, \ h}$ introduced in (\ref{Nal}). Then the map $v\to T_{a, \ h}v=u$ on $H^{2}({\mathbb R})$ defined by problem (\ref{auxnl})
has a unique fixed point $v_{a,h}$. This is the only solution of equation (\ref{id1}) in $H^{2}({\mathbb R})$.

\noindent
b) When $\displaystyle{h=\frac{2\pi n}{\sqrt{a}}, \ n\in {\mathbb Z}, \ n\neq 0}$,  we assume that $xG(x)\in L^{1}({\mathbb R})$,
orthogonality relations (\ref{oc}) hold and $2\sqrt{\pi}N_{a, \ h}l<1$. Then the map $v\to T_{a, \ h}v=u$ on $H^{2}({\mathbb R})$ defined by equation (\ref{auxnl}) possesses a unique fixed point $v_{a,h}$. This is the only solution of problem (\ref{id1}) in $H^{2}({\mathbb R})$.

\noindent
In both cases a) and b) the fixed point $v_{a,h}, \ a>0, \ h\in {\mathbb R}, \ h\neq 0$  is
nontrivial on the real line provided the intersection of supports of the Fourier transforms of functions
$supp \widehat{F(0, x)}\cap supp {\widehat G}$ is a set of
nonzero Lebesgue measure in ${\mathbb R}$.}

\bigskip

Let us introduce the sequence of approximate equations related to problem (\ref{id1}), namely
\begin{equation}
\label{id1m}
\frac{d^{2}u_{m}(x)}{dx^{2}}+au_{m}(x-h)+\int_{-\infty}^{\infty}G_{m}(x-y)F(u_{m}(y),y)dy=0, \quad x\in {\mathbb R}, \quad m\in {\mathbb N},
\end{equation}
where $a>0, \ h\in {\mathbb R}, \ h\neq 0$ are the constants. The sequence of kernels  $\{G_{m}(x)\}_{m=1}^{\infty}$ converges
to $G(x)$  in the appropriate function spaces listed further down. We establish that, under the reasonable auxiliary assumptions,
each of problems (\ref{id1m}) admits  a unique solution $u_{m}(x)\in H^{2}({\mathbb R})$, limiting equation (\ref{id1}) possesses
a unique solution $u(x)\in H^{2}({\mathbb R})$ and $u_{m}(x)\to u(x)$ in $H^{2}({\mathbb R})$ as $m\to \infty$.
 In this case, the solvability relations can be formulated for the iterated  kernels $G_{m}$.
They give us the convergence of the kernels in terms of the Fourier images
 (see the Appendix) and as a consequence the convergence of the solutions (Theorem 5). Our final main statement is as follows.

\bigskip

\noindent
{\bf Theorem 5.}  {\it  Let $m\in {\mathbb N}, \ G_{m}(x): {\mathbb R}\to {\mathbb R}, \ G_{m}(x)\in L^{1}({\mathbb R})$, so that
$G_{m}(x)\to G(x)$ in  $L^{1}({\mathbb R})$ as $m\to \infty$,  the constants  $a>0, \ h\in {\mathbb R}, \ h\neq 0$
and Assumption 3  holds.

\noindent
 a)  If $\displaystyle{h\neq \frac{2\pi n}{\sqrt{a}}, \ n\in {\mathbb Z}}$, assume that
\begin{equation}
\label{2rpnahmle}
2\sqrt{\pi}N_{a, \ h, \ m}l\leq 1-\varepsilon
\end{equation}
for all $m\in {\mathbb N}$ with some fixed $0<\varepsilon<1$ and $N_{a, \ h, \ m}$ defined in (\ref{Na}). Then each equation
(\ref{id1m}) possesses a unique solution $u_{m}(x)\in H^{2}({\mathbb R})$, and limiting problem (\ref{id1}) has
a unique solution $u(x)\in H^{2}({\mathbb R})$.

\noindent
b) When $\displaystyle{h=\frac{2\pi n}{\sqrt{a}}, \ n\in {\mathbb Z}, \ n\neq 0}$,  we assume that
$xG_{m}(x)\in L^{1}({\mathbb R}), \ xG_{m}(x)\to xG(x)$ in $L^{1}({\mathbb R})$ as $m\to \infty$,
orthogonality conditions (\ref{ocm}) are valid and inequality (\ref{2rpnahmle}) holds
for all $m\in {\mathbb N}$ with a certain fixed $0<\varepsilon<1$ as well. Then each problem (\ref{id1m}) admits
a unique solution $u_{m}(x)\in H^{2}({\mathbb R})$, and limiting equation (\ref{id1}) has a unique solution $u(x)\in H^{2}({\mathbb R})$.

\noindent
In both cases a) and b), we have $u_{m}(x)\to u(x)$ in $H^{2}({\mathbb R})$ as $m\to \infty$.

\noindent
The unique solution $u_{m}(x)$  of each equation  (\ref{id1m}) does not vanish identically in ${\mathbb R}$ provided that the intersection of supports of the Fourier images of functions $supp \widehat{F(0, x)}\cap supp {\widehat G_{m}}$ is a set of
nonzero Lebesgue measure on the real line. Similarly, the unique solution $u(x)$  of limiting problem  (\ref{id1}) is nontrivial if the intersection of supports of the Fourier transforms of functions $supp \widehat{F(0, x)}\cap supp {\widehat G}$ is a set of
nonzero Lebesgue measure in ${\mathbb R}$.}


\setcounter{equation}{0}

\section{Solvability of the Differential Equation}

\bigskip

{\it Proof of Theorem 1.} It can be easily verified that it is sufficient to
solve our equation (\ref{eqsh1}) in $L^{2}({\mathbb R})$. Indeed, if $u(x)$
is a square integrable solution of (\ref{eqsh1}) on the real line, we have
$u(x-h)\in L^{2}({\mathbb R})$ as well via (\ref{uxhpn}). Since the right
side of (\ref{eqsh1}) is square integrable as assumed, we obtain
$\displaystyle{\frac{d^{2}u(x)}{dx^{2}}\in L^{2}({\mathbb R})}$, such that
$u(x)\in H^{2}({\mathbb R})$.

Let us suppose that equation (\ref{eqsh1}) admits two solutions
$u_{1}(x), \ u_{2}(x)\in H^{2}({\mathbb R})$. Evidently, their difference
$w(x)=u_{1}(x)-u_{2}(x)\in H^{2}({\mathbb R})$ solves the homogeneous problem
$$
-\frac{d^{2}w(x)}{dx^{2}}-aw(x-h)=0.
$$
Since the operator $L_{h}: H^{2}({\mathbb R})\to L^{2}({\mathbb R})$ (see (\ref{Lh})) does not
have any nontrivial zero modes, $w(x)\equiv 0$ on ${\mathbb R}$.

We apply the standard Fourier transform (\ref{f}) to both sides of equation
(\ref{eqsh1}) and arrive at
\begin{equation}
\label{uhpf}
{\widehat u}(p)=\frac{{\widehat f}(p)}{p^{2}-a \cos(ph)+ia \sin(ph)}.
\end{equation}
For the norm we easily obtain
\begin{equation}
\label{un}
\|u\|_{L^{2}({\mathbb R})}^{2}=\int_{-\infty}^{\infty}\frac{|{\widehat f}(p)|^{2}}
{(p^{2}-a)^{2}+2ap^{2}(1- \cos(ph))}dp.
\end{equation}
Let us first consider the situation when the argument shift constant
$\displaystyle{h\neq \frac{2\pi n}{\sqrt{a}}, \ n\in {\mathbb Z}}$. By means
of (\ref{un}) along with inequality (\ref{alfa}) we easily derive
$$
\|u(x)\|_{L^{2}({\mathbb R})}^{2}\leq \frac{1}{\alpha}\|f(x)\|_{L^{2}({\mathbb R})}^{2}<
\infty
$$
as assumed. Then according to the argument above the unique solution of
our equation $u(x)\in H^{2}({\mathbb R})$, which completes the proof of
the part a) of the theorem.

We turn our attention to the case when
$\displaystyle{h= \frac{2\pi n}{\sqrt{a}}, \ n\in {\mathbb Z}, \ n\neq 0}$.
For the technical purposes we introduce the intervals
\begin{equation}
\label{id}
I_{\delta}^{+}:=[\sqrt{a}-\delta, \ \sqrt{a}+\delta], \quad
I_{\delta}^{-}:=[-\sqrt{a}-\delta, \ -\sqrt{a}+\delta], \quad
0<\delta<\sqrt{a},
\end{equation}
such that their union $I_{\delta}:=I_{\delta}^{+}\cup I_{\delta}^{-}$ and
$I_{\delta}^{c}$ is the complement of the set $I_{\delta}$
on the real line. Let us also use
\begin{equation}
\label{idpm}
I_{\delta}^{c +}:=I_{\delta}^{c}\cap {\mathbb R}^{+},  \quad
I_{\delta}^{c -}:=I_{\delta}^{c}\cap {\mathbb R}^{-}.
\end{equation}
We express
$$
{\widehat u}(p)=\frac{{\widehat f}(p)}{p^{2}-ae^{-iph}}\chi_{I_{\delta}^{+}}+
\frac{{\widehat f}(p)}{p^{2}-ae^{-iph}}\chi_{I_{\delta}^{-}}+
$$
\begin{equation}
\label{uhpf4}
\frac{{\widehat f}(p)}{p^{2}-ae^{-iph}}\chi_{I_{\delta}^{c +}}+
\frac{{\widehat f}(p)}{p^{2}-ae^{-iph}}\chi_{I_{\delta}^{c -}}.
\end{equation}
Here and below $\chi_{A}$ will stand for the characteristic function of a set $A\subseteq {\mathbb R}$.
The third term in the right side of (\ref{uhpf4}) can be estimated from above in the norm as
$$
\Bigg\|\frac{{\widehat f}(p)}{p^{2}-ae^{-iph}}\chi_{I_{\delta}^{c +}}\Bigg\|_{L^{2}({\mathbb R})}^{2}=\int_{-\infty}^{\infty}\frac
{|{\widehat f}(p)|^{2}}{(p^{2}-a)^{2}+2ap^{2}(1-\cos(ph))}\chi_{I_{\delta}^{c +}}dp\leq
$$
$$
\int_{-\infty}^{\infty}\frac{|{\widehat f}(p)|^{2}}{(p-\sqrt{a})^{2}(p+\sqrt{a})^{2}}\chi_{I_{\delta}^{c +}}dp\leq
\frac{1}{{\delta}^{2}a}\|f(x)\|_{L^{2}({\mathbb R})}^{2}<\infty
$$
as we assume, such that
$$
\frac{{\widehat f}(p)}{p^{2}-ae^{-iph}}\chi_{I_{\delta}^{c +}}\in L^{2}({\mathbb R}).
$$
The fourth term in the right side of (\ref{uhpf4}) can be bounded from above in the norm as
$$
\Bigg\|\frac{{\widehat f}(p)}{p^{2}-ae^{-iph}}\chi_{I_{\delta}^{c -}}\Bigg\|_{L^{2}({\mathbb R})}^{2}=\int_{-\infty}^{\infty}\frac
{|{\widehat f}(p)|^{2}}{(p^{2}-a)^{2}+2ap^{2}(1- \cos(ph))}\chi_{I_{\delta}^{c -}}dp\leq
$$
$$
\int_{-\infty}^{\infty}\frac{|{\widehat f}(p)|^{2}}{(p-\sqrt{a})^{2}(p+\sqrt{a})^{2}}\chi_{I_{\delta}^{c -}}dp\leq
\frac{1}{{\delta}^{2}a}\|f(x)\|_{L^{2}({\mathbb R})}^{2}<\infty
$$
as above, so that
$$
\frac{{\widehat f}(p)}{p^{2}-ae^{-iph}}\chi_{I_{\delta}^{c -}}\in L^{2}({\mathbb R}).
$$
Clearly, we have
$$
{\widehat f}(p)={\widehat f}(\sqrt{a})+\int_{\sqrt{a}}^{p}\frac{d{\widehat f}(s)}{ds}ds.
$$
Then the first term in the right side of (\ref{uhpf4}) can be written as
\begin{equation}
\label{fhraidp}
\frac{{\widehat f}(\sqrt{a})}{p^{2}-ae^{-iph}}\chi_{I_{\delta}^{+}}+\frac{\int_{\sqrt{a}}^{p}\frac{d{\widehat f}(s)}{ds}ds}
{p^{2}-ae^{-iph}}\chi_{I_{\delta}^{+}}.
\end{equation}
It can be easily derived from the definition of the standard Fourier transform (\ref{f}) that
\begin{equation}
\label{dfhpdp}
\Bigg|\frac{d{\widehat f}(p)}{dp}\Bigg|\leq \frac{1}{\sqrt{2\pi}}\|xf(x)\|_{L^{1}({\mathbb R})}.
\end{equation}
By means of  (\ref{dfhpdp}), we obtain the estimate from above in the absolute value on the second term in  (\ref{fhraidp}), namely
$$
\Bigg|\frac{\int_{\sqrt{a}}^{p}\frac{d{\widehat f}(s)}{ds}ds}{p^{2}-ae^{-iph}}\chi_{I_{\delta}^{+}}\Bigg|\leq \frac
{\|xf(x)\|_{L^{1}({\mathbb R})}|p-\sqrt{a}|}{\sqrt{2\pi}|p^{2}-ae^{-iph}|}\chi_{I_{\delta}^{+}}.
$$
Then
$$
\Bigg|\frac{\int_{\sqrt{a}}^{p}\frac{d{\widehat f}(s)}{ds}ds}{p^{2}-ae^{-iph}}\chi_{I_{\delta}^{+}}\Bigg|^{2}\leq \frac{1}{2\pi}
\|xf(x)\|_{L^{1}({\mathbb R})}^{2}\frac{|p-\sqrt{a}|^{2}\chi_{I_{\delta}^{+}}}{(p^{2}-a)^{2}+2ap^{2}(1-\cos(ph))}\leq
$$
$$
\frac{1}{2\pi}\|xf(x)\|_{L^{1}({\mathbb R})}^{2}\frac{\chi_{I_{\delta}^{+}}}{(p+\sqrt{a})^{2}}\leq
\frac{1}{2\pi}\|xf(x)\|_{L^{1}({\mathbb R})}^{2}\frac{\chi_{I_{\delta}^{+}}}{a}.
$$
This means that
$$
\frac{\int_{\sqrt{a}}^{p}\frac{d{\widehat f}(s)}{ds}ds}{p^{2}-ae^{-iph}}\chi_{I_{\delta}^{+}}\in L^{2}({\mathbb R})
$$
and it remains to consider the term
\begin{equation}
\label{fhraidpr}
\frac{{\widehat f}(\sqrt{a})}{p^{2}-ae^{-iph}}\chi_{I_{\delta}^{+}}
\end{equation}
in sum (\ref{fhraidp}). Clearly,
$$
\Bigg\|\frac{{\widehat f}(\sqrt{a})}{p^{2}-ae^{-iph}}\chi_{I_{\delta}^{+}}\Bigg\|_{L^{2}({\mathbb R})}^{2}=\int_{-\infty}^{\infty}
\frac{|{\widehat f}(\sqrt{a})|^{2}\chi_{I_{\delta}^{+}}}{(p^{2}-a)^{2}+2ap^{2}\Big(1-\cos\Big(p\frac{2\pi n}{\sqrt{a}}\Big)\Big)}dp=
$$
\begin{equation}
\label{irapmd}
|{\widehat f}(\sqrt{a})|^{2}\int_{\sqrt{a}-\delta}^{\sqrt{a}+\delta}\frac{dp}{(p-\sqrt{a})^{2}(p+\sqrt{a})^{2}+
4ap^{2}\sin^{2}\Big(\frac{\pi np}{\sqrt{a}}\Big)}.
\end{equation}
By virtue of the change of variables
\begin{equation}
\label{xi}
p=\sqrt{a}\xi,
\end{equation}
the right side of  (\ref{irapmd}) equals to
\begin{equation}
\label{irapmdxi}
\frac{|{\widehat f}(\sqrt{a})|^{2}}{a^{\frac{3}{2}}}\int_{\frac{\sqrt{a}-\delta}{\sqrt{a}}}^{\frac{\sqrt{a}+\delta}{\sqrt{a}}}\frac{d\xi}
{(\xi-1)^{2}(\xi+1)^{2}+4\xi^{2}sin^{2}(\pi n[\xi-1])}.
\end{equation}
Let $\eta:=\xi-1$. Then (\ref{irapmdxi}) is equal to
\begin{equation}
\label{irapmdeta}
\frac{|{\widehat f}(\sqrt{a})|^{2}}{a^{\frac{3}{2}}}\int_{-\frac{\delta}{\sqrt{a}}}^{\frac{\delta}{\sqrt{a}}}\frac{d\eta}
{\eta^{2}(\eta+2)^{2}+4(\eta+1)^{2}\sin^{2}(\pi n\eta)}.
\end{equation}
Note that for
$\displaystyle{\eta\in \Big[-\frac{\delta}{\sqrt{a}}, \frac{\delta}{\sqrt{a}}\Big]\subset [-1, 1]}$, we have $(\eta+2)^{2}\leq 9$ and
$(\eta+1)^{2}\leq 4$, such that (\ref{irapmdeta}) can be easily bounded from below by
\begin{equation}
\label{irapmdetan+}
\frac{2|{\widehat f}(\sqrt{a})|^{2}}{a^{\frac{3}{2}}}\int_{0}^{\frac{\delta}{\sqrt{a}}}\frac{d\eta}{9\eta^{2}+16\sin^{2}(\pi n\eta)}.
\end{equation}
Evidently,
\begin{equation}
\label{sin}
sin^{2}(\pi n\eta)\leq \pi^{2}n^{2}\eta^{2}, \ \eta\in \Big[0, \frac{\delta}{\sqrt{a}}\Big].
\end{equation}
Therefore, expression (\ref{irapmdetan+}) is
infinite unless ${\widehat f}(\sqrt{a})$ vanishes. This is equivalent to the orthogonality condition
\begin{equation}
\label{oc1}
\Big(f(x), \frac{e^{i\sqrt{a}x}}{\sqrt{2\pi}}\Big)_{L^{2}({\mathbb R})}=0.
\end{equation}
Obviously,
$$
{\widehat f}(p)={\widehat f}(-\sqrt{a})+\int_{-\sqrt{a}}^{p}\frac{d{\widehat f}(s)}{ds}ds.
$$
Hence, the second term in the right side of  (\ref{uhpf4}) is given by
\begin{equation}
\label{fhmra2i}
\frac{{\widehat f}(-\sqrt{a})}{p^{2}-ae^{-iph}}\chi_{I_{\delta}^{-}}+
\frac{\int_{-\sqrt{a}}^{p}\frac{d{\widehat f}(s)}{ds}ds}{p^{2}-ae^{-iph}}\chi_{I_{\delta}^{-}}.
\end{equation}
Let us recall inequality (\ref{dfhpdp}). Thus,
$$
\Bigg|\frac{\int_{-\sqrt{a}}^{p}\frac{d{\widehat f}(s)}{ds}ds}{p^{2}-ae^{-iph}}\chi_{I_{\delta}^{-}}\Bigg|\leq
\frac{\|xf(x)\|_{L^{1}({\mathbb R})}|p+\sqrt{a}|}{\sqrt{2\pi}|p^{2}-ae^{-iph}|}\chi_{I_{\delta}^{-}},
$$
so that
$$
\Bigg|\frac{\int_{-\sqrt{a}}^{p}\frac{d{\widehat f}(s)}{ds}ds}{p^{2}-ae^{-iph}}\chi_{I_{\delta}^{-}}\Bigg|^{2}\leq \frac{1}{2\pi}
\|xf(x)\|_{L^{1}({\mathbb R})}^{2}\frac{|p+\sqrt{a}|^{2}\chi_{I_{\delta}^{-}}}
{(p^{2}-a)^{2}+2ap^{2}(1-\cos(ph))}\leq
$$
$$
\frac{\|xf(x)\|_{L^{1}({\mathbb R})}^{2}}{2\pi}\frac{\chi_{I_{\delta}^{-}}}{(p-\sqrt{a})^{2}}\leq
\frac{\|xf(x)\|_{L^{1}({\mathbb R})}^{2}}{2\pi}\frac{\chi_{I_{\delta}^{-}}}{a}.
$$
Therefore,
$$
\frac{\int_{-\sqrt{a}}^{p}\frac{d{\widehat f}(s)}{ds}ds}{p^{2}-ae^{-iph}}\chi_{I_{\delta}^{-}}\in L^{2}({\mathbb R}).
$$
Finally, we analyze the term
\begin{equation}
\label{fhmradm}
\frac{{\widehat f}(-\sqrt{a})}{p^{2}-ae^{-iph}}\chi_{I_{\delta}^{-}}
\end{equation}
in sum (\ref{fhmra2i}). We have
$$
\Bigg\|\frac{{\widehat f}(-\sqrt{a})}{p^{2}-ae^{-iph}}\chi_{I_{\delta}^{-}}\Bigg\|_{L^{2}({\mathbb R})}^{2}=
\int_{-\infty}^{\infty}
\frac{|{\widehat f}(-\sqrt{a})|^{2}\chi_{I_{\delta}^{-}}}{(p^{2}-a)^{2}+2ap^{2}\Big(1-\cos\Big(p\frac{2\pi n}{\sqrt{a}}\Big)\Big)}dp=
$$
\begin{equation}
\label{irapmdn}
|{\widehat f}(-\sqrt{a})|^{2}\int_{-\sqrt{a}-\delta}^{-\sqrt{a}+\delta}\frac{dp}{(p-\sqrt{a})^{2}(p+\sqrt{a})^{2}+
4ap^{2}\sin^{2}\Big(\frac{\pi np}{\sqrt{a}}\Big)}.
\end{equation}
Recall the change of variables (\ref{xi}). Hence, the right side of (\ref{irapmdn}) is given by
\begin{equation}
\label{irapmdxin}
\frac{|{\widehat f}(-\sqrt{a})|^{2}}{a^{\frac{3}{2}}}\int_{\frac{-\sqrt{a}-\delta}{\sqrt{a}}}^{\frac{-\sqrt{a}+\delta}{\sqrt{a}}}\frac{d\xi}
{(\xi-1)^{2}(\xi+1)^{2}+4\xi^{2}\sin^{2}(\pi n[\xi+1])}.
\end{equation}
Let $\gamma:=\xi+1$. Thus, (\ref{irapmdxin}) equals to
\begin{equation}
\label{irapmdetan}
\frac{|{\widehat f}(-\sqrt{a})|^{2}}{a^{\frac{3}{2}}}\int_{-\frac{\delta}{\sqrt{a}}}^{\frac{\delta}{\sqrt{a}}}\frac{d\gamma}
{\gamma^{2}(\gamma-2)^{2}+4(\gamma-1)^{2}\sin^{2}(\pi n\gamma)}.
\end{equation}
Clearly,  for
$\displaystyle{\gamma\in \Big[-\frac{\delta}{\sqrt{a}}, \frac{\delta}{\sqrt{a}}\Big]\subset [-1, 1]}$, we obtain $(\gamma-2)^{2}\leq 9$ and
$(\gamma-1)^{2}\leq 4$, so that (\ref{irapmdetan}) can be easily estimated from below by
\begin{equation}
\label{irapmdetann}
\frac{2|{\widehat f}(-\sqrt{a})|^{2}}{a^{\frac{3}{2}}}\int_{0}^{\frac{\delta}{\sqrt{a}}}\frac{d\gamma}{9\gamma^{2}+16\sin^{2}(\pi n\gamma)}.
\end{equation}
Let us use here the trivial bound analogous to (\ref{sin}). Thus, (\ref{irapmdetann}) is not finite unless
${\widehat f}(-\sqrt{a})=0$, which is equivalent to the orthogonality relation
\begin{equation}
\label{oc2}
\Big(f(x), \frac{e^{-i\sqrt{a}x}}{\sqrt{2\pi}}\Big)_{L^{2}({\mathbb R})}=0.
\end{equation}
Formulas (\ref{oc1}) and (\ref{oc2}) give the solvability conditions for our equation  (\ref{eqsh1}) in the situation when
$\displaystyle{h= \frac{2\pi n}{\sqrt{a}}, \ n\in {\mathbb Z}, \ n\neq 0}$.  Note that the terms (\ref{fhraidpr}) and (\ref{fhmradm}) are
orthogonal to each other in $L^{2}({\mathbb R})$ since they have disjoint supports. \hfill\lanbox

\bigskip

Then we turn our attention to establishing the solvability in the sense of sequences for our problem (\ref{eqsh1}).

\bigskip

\noindent
{\it Proof of Theorem 2.} It can be trivially checked that if each equation (\ref{eqsh1m}) admits a unique solution
$u_{m}(x)\in H^{2}({\mathbb R})$, limiting problem (\ref{eqsh1}) possesses a unique solution $u(x)\in H^{2}({\mathbb R})$ and it is known
that $u_{m}(x)\to u(x)$ in $L^{2}({\mathbb R})$ as $m\to \infty$, then we have
$u_{m}(x)\to u(x)$ in $H^{2}({\mathbb R})$ as $m\to \infty$. Indeed, from formulas (\ref{eqsh1m})  and (\ref{eqsh1}) we easily deduce
that
$$
-\Big(\frac{d^{2}u_{m}(x)}{dx^{2}}-\frac{d^{2}u(x)}{dx^{2}}\Big)=a[u_{m}(x-h)-u(x-h)]+f_{m}(x)-f(x),
$$
such that
\begin{equation}
\label{umul2r}
\Big\|\frac{d^{2}u_{m}(x)}{dx^{2}}-\frac{d^{2}u(x)}{dx^{2}}\Big\|_{L^{2}({\mathbb R})}\leq a\|u_{m}(x)-u(x)\|_{L^{2}({\mathbb R})}+
\|f_{m}(x)-f(x)\|_{L^{2}({\mathbb R})}.
\end{equation}
Here we used the translation invariance of the $L^{2}({\mathbb R})$ norm. Then according to our assumptions we have
$$
\frac{d^{2}u_{m}(x)}{dx^{2}} \to \frac{d^{2}u(x)}{dx^{2}} \quad in \quad L^{2}({\mathbb R})
$$
as $m\to \infty$. Let us recall the definition of the norm (\ref{sob}). Therefore,
$$
u_{m}(x)\to u(x) \quad in \quad H^{2}({\mathbb R})
$$
as $m\to \infty$ as well.

First we deal with the case when the argument translation constant
$\displaystyle{h\neq \frac{2\pi n}{\sqrt{a}}, \ n\in {\mathbb Z}}$.  By virtue of the result of Theorem 1 above,  equations
(\ref{eqsh1m})  and (\ref{eqsh1}) admit unique solutions $u_{m}(x)\in H^{2}({\mathbb R}), \ m\in {\mathbb N}$ and
$u(x)\in H^{2}({\mathbb R})$ respectively. Let us apply the standard Fourier transform (\ref{f}) to both sides of problems
(\ref{eqsh1m})  and (\ref{eqsh1}). This yields
\begin{equation}
\label{umuhp}
{\widehat u}(p)=\frac{{\widehat f}(p)}{p^{2}-ae^{-iph}}, \quad
{\widehat u_{m}}(p)=\frac{{\widehat f_{m}}(p)}{p^{2}-ae^{-iph}}, \quad m\in {\mathbb N},
\end{equation}
so that
\begin{equation}
\label{umuhpn}
\|u_{m}-u\|_{L^{2}({\mathbb R})}^{2}=\int_{-\infty}^{\infty}\frac{|{\widehat f_{m}}(p)-{\widehat f}(p)|^{2}}
{(p^{2}-a)^{2}+2ap^{2}(1-\cos(ph))}dp.
\end{equation}
Let us recall estimate (\ref{alfa}). By virtue of  (\ref{umuhpn}), we have
$$
\|u_{m}(x)-u(x)\|_{L^{2}({\mathbb R})}^{2}\leq \frac{1}{\alpha}\|f_{m}(x)-f(x)\|_{L^{2}({\mathbb R})}^{2}\to 0
$$
as $m\to \infty$ as we assume. Hence,
$$
u_{m}(x)\to u(x) \quad in \quad L^{2}({\mathbb R})
$$
as $m\to \infty$. By means of the reasoning above,
$$
u_{m}(x)\to u(x) \quad in \quad H^{2}({\mathbb R})
$$
as $m\to \infty$, which completes the proof of our theorem in the situation a).

Then we consider the case when $\displaystyle{h=\frac{2\pi n}{\sqrt{a}}, \ n\in {\mathbb Z}, \ n\neq 0}$. According to the result of
Theorem 1, each equation (\ref{eqsh1m}) has a unique solution $u_{m}(x)\in H^{2}({\mathbb R})$. Recall the part a) of Lemma 3.3 of  ~\cite{VV1031}. Therefore, under our assumptions, the orthogonality conditions
\begin{equation}
\label{ocl}
\Big(f(x), \frac{e^{\pm i \sqrt{a}x}}{\sqrt{2\pi}}\Big)_{L^{2}({\mathbb R})}=0
\end{equation}
hold. By virtue of Theorem 1, limiting problem (\ref{eqsh1}) admits a unique solution $u(x)\in H^{2}({\mathbb R})$. We apply the
standard Fourier transform to both sides of equations (\ref{eqsh1m}) and (\ref{eqsh1})  and derive  (\ref{umuhp}). Note that
$$
 {\widehat u_{m}}(p) -{\widehat u}(p)=\frac{{\widehat f_{m}}(p)-{\widehat f}(p)}{p^{2}-ae^{-iph}}\chi_{I_{\delta}^{+}}+
\frac{{\widehat f_{m}}(p) -{\widehat f}(p)}{p^{2}-ae^{-iph}}\chi_{I_{\delta}^{-}}+
$$
\begin{equation}
\label{uhpf4n}
\frac{{\widehat f_{m}}(p)-{\widehat f}(p)}{p^{2}-ae^{-iph}}\chi_{I_{\delta}^{c +}}+
\frac{ {\widehat f}_{m}(p)-{\widehat f}(p)}{p^{2}-ae^{-iph}}\chi_{I_{\delta}^{c -}}.
\end{equation}
The third term in the right side of  (\ref{uhpf4n}) can be bounded from above in the norm as
$$
\Bigg\|\frac{{\widehat f_{m}}(p)-{\widehat f}(p)}{p^{2}-ae^{-iph}}\chi_{I_{\delta}^{c +}}\Bigg\|_{L^{2}({\mathbb R})}^{2}=\int_{-\infty}^{\infty}
\frac{|{\widehat f_{m}}(p)-{\widehat f}(p)|^{2}}{(p^{2}-a)^{2}+2ap^{2}(1-\cos(ph))}\chi_{I_{\delta}^{c +}}dp\leq
$$
$$
\int_{-\infty}^{\infty}\frac{|{\widehat f_{m}}(p)-{\widehat f}(p)|^{2}}{(p-\sqrt{a})^{2}(p+\sqrt{a})^{2}}\chi_{I_{\delta}^{c +}}dp\leq
\frac{1}{{\delta}^{2}a}\|f_{m}(x)-f(x)\|_{L^{2}({\mathbb R})}^{2}\to 0
$$
as $m\to \infty$ according to our assumption.

The fourth term in the right side of  (\ref{uhpf4n}) can be estimated from above in the norm as
$$
\Bigg\|\frac{{\widehat f_{m}}(p)-{\widehat f}(p)}{p^{2}-ae^{-iph}}\chi_{I_{\delta}^{c -}}\Bigg\|_{L^{2}({\mathbb R})}^{2}=\int_{-\infty}^{\infty}
\frac{|{\widehat f_{m}}(p)-{\widehat f}(p)|^{2}}{(p^{2}-a)^{2}+2ap^{2}(1-\cos(ph))}\chi_{I_{\delta}^{c -}}dp\leq
$$
$$
\int_{-\infty}^{\infty}\frac{|{\widehat f_{m}}(p)-{\widehat f}(p)|^{2}}{(p-\sqrt{a})^{2}(p+\sqrt{a})^{2}}\chi_{I_{\delta}^{c -}}dp\leq
\frac{1}{{\delta}^{2}a}\|f_{m}(x)-f(x)\|_{L^{2}({\mathbb R})}^{2}\to 0
$$
as $m\to \infty$ analogously. Let us recall orthogonality relations (\ref{ocfm})  and (\ref{ocl}). Thus,
$$
{\widehat f}(\sqrt{a})=0, \quad {\widehat f_{m}}(\sqrt{a})=0, \quad m\in {\mathbb N},
$$
such that
$$
{\widehat f}(p)=\int_{\sqrt{a}}^{p}\frac{d{\widehat f}(s)}{ds}ds, \quad
{\widehat f_{m}}(p)=\int_{\sqrt{a}}^{p}\frac{d{\widehat f_{m}}(s)}{ds}ds, \quad m\in {\mathbb N}.
$$
Therefore, the first term in the right side of  (\ref{uhpf4n}) is given by
\begin{equation}
\label{fmhfhi}
\frac{\int_{\sqrt{a}}^{p}\Big[\frac{d{\widehat f_{m}}(s)}{ds}-\frac{d{\widehat f}(s)}{ds}\Big]ds}{p^{2}-ae^{-iph}}\chi_{I_{\delta}^{+}}.
\end{equation}
We use the definition of the standard Fourier transform (\ref{f}) to obtain that
\begin{equation}
\label{dfmhfhm}
\Bigg|\frac{d{\widehat f_{m}}(p)}{dp}-\frac{d{\widehat f}(p)}{dp}\Bigg|\leq \frac{1}{\sqrt{2\pi}}\|xf_{m}(x)-xf(x)\|_{L^{1}({\mathbb R})}.
\end{equation}
By virtue of (\ref{dfmhfhm}), we derive the upper bound on expression (\ref{fmhfhi}) in the absolute value as
$$
\Bigg|\frac{\int_{\sqrt{a}}^{p}\Big[\frac{d{\widehat f_{m}}(s)}{ds}-\frac{d{\widehat f}(s)}{ds}\Big]ds}{p^{2}-ae^{-iph}}\chi_{I_{\delta}^{+}}\Bigg|
\leq \frac{\|xf_{m}(x)-xf(x)\|_{L^{1}({\mathbb R})}|p-\sqrt{a}|}{\sqrt{2\pi}|p^{2}-ae^{-iph}|}\chi_{I_{\delta}^{+}}.
$$
Clearly,
$$
\Bigg|\frac{\int_{\sqrt{a}}^{p}\Big[\frac{d{\widehat f_{m}}(s)}{ds}-\frac{d{\widehat f}(s)}{ds}\Big]ds}{p^{2}-ae^{-iph}}\chi_{I_{\delta}^{+}}\Bigg|^{2}\leq \frac{1}{2\pi}\|xf_{m}(x)-xf(x)\|_{L^{1}({\mathbb R})}^{2}\frac{|p-\sqrt{a}|^{2}\chi_{I_{\delta}^{+}}}
{(p^{2}-a)^{2}+2ap^{2}(1-\cos(ph))}\leq
$$
$$
\frac{1}{2\pi}\frac{\|xf_{m}(x)-xf(x)\|_{L^{1}({\mathbb R})}^{2}\chi_{I_{\delta}^{+}}}{(p+\sqrt{a})^{2}}\leq
\frac{1}{2\pi}\frac{\|xf_{m}(x)-xf(x)\|_{L^{1}({\mathbb R})}^{2}\chi_{I_{\delta}^{+}}}{a},
$$
so that
$$
\Bigg\|\frac{\int_{\sqrt{a}}^{p}\Big[\frac{d{\widehat f_{m}}(s)}{ds}-\frac{d{\widehat f}(s)}{ds}\Big]ds}{p^{2}-ae^{-iph}}\chi_{I_{\delta}^{+}}
\Bigg\|_{L^{2}({\mathbb R})}\leq \sqrt{\frac{\delta}{\pi a}}\|xf_{m}(x)-xf(x)\|_{L^{1}({\mathbb R})}\to 0
$$
as $m\to \infty$ due to the given condition. By means of orthogonality relations (\ref{ocfm}) and (\ref{ocl}), we have
$$
{\widehat f}(-\sqrt{a})=0, \quad {\widehat f_{m}}(-\sqrt{a})=0, \quad m\in{\mathbb N},
$$
such that
$$
{\widehat f}(p)=\int_{-\sqrt{a}}^{p}\frac{d{\widehat f}(s)}{ds}ds, \quad
{\widehat f_{m}}(p)=\int_{-\sqrt{a}}^{p}\frac{d{\widehat f_{m}}(s)}{ds}ds, \quad m\in{\mathbb N}.
$$
Hence, the second term in the right side of  (\ref{uhpf4n}) equals to
\begin{equation}
\label{fmhfhim}
\frac{\int_{-\sqrt{a}}^{p}\Big[\frac{d{\widehat f_{m}}(s)}{ds}-\frac{d{\widehat f}(s)}{ds}\Big]ds}{p^{2}-ae^{-iph}}\chi_{I_{\delta}^{-}}.
\end{equation}
By virtue of inequality (\ref{dfmhfhm}), we arrive at the estimate from above
$$
\Bigg|\frac{\int_{-\sqrt{a}}^{p}\Big[\frac{d{\widehat f_{m}}(s)}{ds}-\frac{d{\widehat f}(s)}{ds}\Big]ds}{p^{2}-ae^{-iph}}\chi_{I_{\delta}^{-}}\Bigg|
\leq \frac{\|xf_{m}(x)-xf(x)\|_{L^{1}({\mathbb R})}|p+\sqrt{a}|}{\sqrt{2\pi}|p^{2}-ae^{-iph}|}\chi_{I_{\delta}^{-}}.
$$
Evidently,
$$
\Bigg|\frac{\int_{-\sqrt{a}}^{p}\Big[\frac{d{\widehat f_{m}}(s)}{ds}-\frac{d{\widehat f}(s)}{ds}\Big]ds}{p^{2}-ae^{-iph}}\chi_{I_{\delta}^{-}}\Bigg|^{2}\leq \frac{1}{2\pi}\|xf_{m}(x)-xf(x)\|_{L^{1}({\mathbb R})}^{2}\frac{|p+\sqrt{a}|^{2}\chi_{I_{\delta}^{-}}}
{(p^{2}-a)^{2}+2ap^{2}(1-\cos(ph))}\leq
$$
$$
\frac{1}{2\pi}\frac{\|xf_{m}(x)-xf(x)\|_{L^{1}({\mathbb R})}^{2}\chi_{I_{\delta}^{-}}}{(p-\sqrt{a})^{2}}\leq
\frac{1}{2\pi}\frac{\|xf_{m}(x)-xf(x)\|_{L^{1}({\mathbb R})}^{2}\chi_{I_{\delta}^{-}}}{a},
$$
so that
$$
\Bigg\|\frac{\int_{-\sqrt{a}}^{p}\Big[\frac{d{\widehat f_{m}}(s)}{ds}-\frac{d{\widehat f}(s)}{ds}\Big]ds}{p^{2}-ae^{-iph}}\chi_{I_{\delta}^{-}}
\Bigg\|_{L^{2}({\mathbb R})}\leq \sqrt{\frac{\delta}{\pi a}}\|xf_{m}(x)-xf(x)\|_{L^{1}({\mathbb R})}\to 0
$$
as $m\to \infty$ as above. This means that
$$
u_{m}(x)\to u(x) \quad in \quad L^{2}({\mathbb R})
$$
as $m\to \infty$. According to the argument above, we have
$$
u_{m}(x)\to u(x) \quad in \quad H^{2}({\mathbb R})
$$
as $m\to \infty$ in the situation b) of the theorem as well.  \hfill\lanbox


 \setcounter{equation}{0}

 \section{Solvability of the Integro-Differential Equation}

\bigskip

\noindent
{\it Proof of Theorem 4.} Let us first suppose that for a certain $v(x)\in H^{2}({\mathbb R})$ there exist two solutions
$u_{1, 2}(x)\in H^{2}({\mathbb R})$ of problem (\ref{auxnl}). Obviously, their difference
$w(x):=u_{1}(x)-u_{2}(x)\in H^{2}({\mathbb R})$ will be a solution of the homogeneous equation
$$
-\frac{d^{2}w(x)}{dx^{2}}-aw(x-h)=0.
$$
The operator (\ref{Lh}) involved in its left side does not have any nontrivial zero modes. Therefore,
$w(x)\equiv 0$ on the real line.

We choose an arbitrary $v(x)\in H^{2}({\mathbb R})$, apply the standard Fourier transform (\ref{f})  to both sides of
(\ref{auxnl}) and arrive at
\begin{equation}
\label{uhpghfheph}
\widehat{u}(p)=\sqrt{2\pi}\frac{\widehat{G}(p)\widehat{f}(p)}{p^{2}-ae^{-iph}}, \quad
p^{2}\widehat{u}(p)=\sqrt{2\pi}\frac{p^{2}\widehat{G}(p)\widehat{f}(p)}{p^{2}-ae^{-iph}}.
\end{equation}
Here $\widehat{f}(p)$ denotes the Fourier image of $F(v(x), x)$. Clearly, the estimates from above
$$
|\widehat{u}(p)|\leq \sqrt{2\pi} N_{a, \ h}|\widehat{f}(p)| \quad and \quad
|p^{2}\widehat{u}(p)|\leq \sqrt{2\pi} N_{a, \ h}|\widehat{f}(p)|
$$
are valid. Note that $N_{a, \ h}<\infty$ by means of the result of Lemma A1 of the Appendix without any orthogonality
relations in the situation a) of our theorem and under orthogonality conditions (\ref{oc}) in the case b). This enables us to obtain
the upper bound on the norm
$$
\|u\|_{H^{2}({\mathbb R})}^{2}=\|\widehat{u}(p)\|_{L^{2}({\mathbb R})}^{2}+
\|p^{2}\widehat{u}(p)\|_{L^{2}({\mathbb R})}^{2}\leq
$$
\begin{equation}
\label{uhpp2uhpl2n}
4\pi N_{a, \ h}^{2}\int_{-\infty}^{\infty}|\widehat{f}(p)|^{2}dp=4\pi N_{a, \ h}^{2}\|F(v(x), x)\|_{L^{2}({\mathbb R})}^{2}.
\end{equation}
Recall condition (\ref{ub1}) of our Assumption 3. Thus, the right side of  (\ref{uhpp2uhpl2n}) is finite for
$v(x)\in L^{2}({\mathbb R})$. Therefore, for an arbitrary $v(x)\in H^{2}({\mathbb R})$  there exists a unique solution
 $u(x)\in H^{2}({\mathbb R})$ of problem (\ref{auxnl}). Its Fourier image is given by (\ref{uhpghfheph}). Hence, the map
$T_{a, \ h}: H^{2}({\mathbb R})\to H^{2}({\mathbb R})$ is well defined.

This alows us to choose arbitrarily $v_{1, 2}(x)\in H^{2}({\mathbb R})$, such that their images
$u_{1, 2}:=T_{a, \ h}v_{1, 2}\in H^{2}({\mathbb R})$. According to (\ref{auxnl}), we have
\begin{equation}
\label{auxnl1}
-\frac{d^{2}u_{1}(x)}{dx^{2}}-au_{1}(x-h)=\int_{-\infty}^{\infty}G(x-y)F(v_{1}(y), y)dy, \quad x\in {\mathbb R},
\end{equation}
\begin{equation}
\label{auxnl2}
-\frac{d^{2}u_{2}(x)}{dx^{2}}-au_{2}(x-h)=\int_{-\infty}^{\infty}G(x-y)F(v_{2}(y), y)dy, \quad x\in {\mathbb R},
\end{equation}
where $a>0, \ h\in {\mathbb R}, \ h\neq 0$ are the constants. Let us apply the standard Fourier transform (\ref{f}) to both sides of
equations (\ref{auxnl1}) and (\ref{auxnl2}). We derive
\begin{equation}
\label{uhpghfheph}
\widehat{u_{1}}(p)=\sqrt{2\pi}\frac{\widehat{G}(p)\widehat{f_{1}}(p)}{p^{2}-ae^{-iph}}, \quad
p^{2}\widehat{u_{1}}(p)=\sqrt{2\pi}\frac{p^{2}\widehat{G}(p)\widehat{f_{1}}(p)}{p^{2}-ae^{-iph}}.
\end{equation}
\begin{equation}
\label{uhpghfheph}
\widehat{u_{2}}(p)=\sqrt{2\pi}\frac{\widehat{G}(p)\widehat{f_{2}}(p)}{p^{2}-ae^{-iph}}, \quad
p^{2}\widehat{u_{2}}(p)=\sqrt{2\pi}\frac{p^{2}\widehat{G}(p)\widehat{f_{2}}(p)}{p^{2}-ae^{-iph}}.
\end{equation}
Here $\widehat{f_{1}}(p)$ and  $\widehat{f_{2}}(p)$ stand for the Fourier images of $F(v_{1}(x), x)$ and $F(v_{2}(x), x)$ respectively.
Evidently, the estimates from above
$$
|\widehat{u_{1}}(p)-\widehat{u_{2}}(p)|\leq \sqrt{2\pi}N_{a, \ h}|\widehat{f_{1}}(p)-\widehat{f_{2}}(p)|, \quad
|p^{2}\widehat{u_{1}}(p)-p^{2}\widehat{u_{2}}(p)|\leq \sqrt{2\pi}N_{a, \ h}|\widehat{f_{1}}(p)-\widehat{f_{2}}(p)|
$$
hold. Then we obtain the upper bound on the norm
$$
\|u_{1}-u_{2}\|_{H^{2}({\mathbb R})}^{2}=\|\widehat{u_{1}}(p)-\widehat{u_{2}}(p)\|_{L^{2}({\mathbb R})}^{2}+
\|p^{2}[\widehat{u_{1}}(p)-\widehat{u_{2}}(p)]\|_{L^{2}({\mathbb R})}^{2}\leq
$$
$$
4\pi N_{a, \ h}^{2}\int_{-\infty}^{\infty}|\widehat{f_{1}}(p)-\widehat{f_{2}}(p)|^{2}dp=4\pi N_{a, \ h}^{2}
\|F(v_{1}(x), x)-F(v_{2}(x), x)\|_{L^{2}({\mathbb R})}^{2}.
$$
Clearly, $v_{1, 2}(x)\in H^{2}({\mathbb R})\subset L^{\infty}({\mathbb R})$ via the Sobolev embedding. Let us recall inequality
(\ref{lk1}) of Assumption 3 above. This yields
$$
\|F(v_{1}(x), x)-F(v_{2}(x), x)\|_{L^{2}({\mathbb R})}\leq l\|\\v_{1}(x)-v_{2}(x)\|_{L^{2}({\mathbb R})},
$$
so that
\begin{equation}
\label{contr}
\|T_{a, \ h}v_{1}-T_{a, \  h}v_{2}\|_{H^{2}({\mathbb R})}\leq 2\sqrt{\pi}N_{a, \ h}l\|v_{1}-v_{2}\|_{H^{2}({\mathbb R})}.
\end{equation}
Note that the constant in the righ tside of estimate (\ref{contr})  is less than one according to the given condition.
By means of the Fixed Point Theorem, there exists a unique function $v_{a, h}\in H^{2}({\mathbb R})$ with
the property $T_{a, \ h}v_{a, h}=v_{a, h}$. This is the only solution of equation (\ref{id1}) in $H^{2}({\mathbb R})$. Suppose
$v_{a, h}\equiv 0$ on the real line. This will contradict to our assumption that the Fourier transforms of $G(x)$ and $F(0, x)$ do not
vanish identically on a set of nonzero Lebesgue measure in ${\mathbb R}$.  \hfill\lanbox

\bigskip

Let us demonstrate the validity of the final main proposition of the article.

\bigskip

\noindent
{\it Proof of Theorem 5.} By means of the result of Theorem 4, each problem (\ref{id1m}) has a unique solution
$u_{m}(x)\in H^{2}({\mathbb R}),  \ m\in {\mathbb N}$. Limiting equation (\ref{id1}) possesses a unique solution
$u(x)\in H^{2}({\mathbb R})$ by virtue of Lemma A2 below and Theorem 4. We apply the standard Fourier transform
(\ref{f}) to both sides of  (\ref{id1}) and (\ref{id1m})  and obtain
\begin{equation}
\label{uhpumhpg}
\widehat{u}(p)=\sqrt{2\pi}\frac{\widehat{G}(p)\widehat{g}(p)}{p^{2}-ae^{-iph}}, \quad
\widehat{u_{m}}(p)=\sqrt{2\pi}\frac{\widehat{G_{m}}(p)\widehat{g_{m}}(p)}{p^{2}-ae^{-iph}}, \quad m\in {\mathbb N}.
\end{equation}
Here $\widehat{g}(p)$ and $\widehat{g_{m}}(p)$ stand for the Fourier images of $F(u(x), x)$ and $F(u_{m}(x), x)$ respectively.
Evidently,
$$
|\widehat{u_{m}}(p)-\widehat{u}(p)|\leq \sqrt{2\pi}\Bigg\|\frac{\widehat{G_{m}}(p)}{p^{2}-ae^{-iph}}-
\frac{\widehat{G}(p)}{p^{2}-ae^{-iph}}\Bigg\|_{L^{\infty}({\mathbb R})}|\widehat{g}(p)|+
$$
$$
\sqrt{2\pi}\Bigg\|\frac{\widehat{G_{m}}(p)}{p^{2}-ae^{-iph}}\Bigg\|_{L^{\infty}({\mathbb R})}|\widehat{g_{m}}(p)-\widehat{g}(p)|.
$$
Hence,
$$
\|u_{m}-u\|_{L^{2}({\mathbb R})}\leq \sqrt{2\pi}\Bigg\|\frac{\widehat{G_{m}}(p)}{p^{2}-ae^{-iph}}-
\frac{\widehat{G}(p)}{p^{2}-ae^{-iph}}\Bigg\|_{L^{\infty}({\mathbb R})}\|F(u(x), x)\|_{L^{2}({\mathbb R})}+
$$
$$
\sqrt{2\pi}\Bigg\|\frac{\widehat{G_{m}}(p)}{p^{2}-ae^{-iph}}\Bigg\|_{L^{\infty}({\mathbb R})}
\|F(u_{m}(x), x)-F(u(x), x)\|_{L^{2}({\mathbb R})}.
$$
Let us recall condition (\ref{lk1}) of Assumption 3. This yields
\begin{equation}
\label{fumuxn}
\|F(u_{m}(x), x)-F(u(x), x)\|_{L^{2}({\mathbb R})}\leq l\|u_{m}(x)-u(x)\|_{L^{2}({\mathbb R})}.
\end{equation}
Clearly, $u_{m}(x), u(x)\in H^{2}({\mathbb R})\subset L^{\infty}({\mathbb R})$ due to the Sobolev embedding. Then, we arrive at
$$
\|u_{m}(x)-u(x)\|_{L^{2}({\mathbb R})}\Bigg\{1-\sqrt{2\pi}\Bigg\|\frac{\widehat{G_{m}}(p)}{p^{2}-ae^{-iph}}\Bigg\|_{L^{\infty}({\mathbb R})}l
\Bigg\}\leq
$$
$$
\sqrt{2\pi}\Bigg\|\frac{\widehat{G_{m}}(p)}{p^{2}-ae^{-iph}}-\frac{\widehat{G}(p)}{p^{2}-ae^{-iph}}\Bigg\|_{L^{\infty}({\mathbb R})}
\|F(u(x), x)\|_{L^{2}({\mathbb R})}.
$$
Let us use inequality (\ref{2rpnahmle}). Thus,
$$
\|u_{m}(x)-u(x)\|_{L^{2}({\mathbb R})}\leq \frac{\sqrt{2\pi}}{\varepsilon}
\Bigg\|\frac{\widehat{G_{m}}(p)}{p^{2}-ae^{-iph}}-\frac{\widehat{G}(p)}{p^{2}-ae^{-iph}}\Bigg\|_{L^{\infty}({\mathbb R})}
\|F(u(x), x)\|_{L^{2}({\mathbb R})}.
$$
By virtue of condition (\ref{ub1}) of Assumption 3, we have $F(u(x), x)\in L^{2}({\mathbb R})$ for
$u(x)\in H^{2}({\mathbb R})$. Therefore,
\begin{equation}
\label{umtoul2}
u_{m}(x)\to u(x), \quad m\to \infty
\end{equation}
in $L^{2}({\mathbb R})$ via the result of Lemma A2 below. Recall  (\ref{uhpumhpg}). Hence,
\begin{equation}
\label{uhpumhpg2}
p^{2}\widehat{u}(p)=\sqrt{2\pi}\frac{p^{2}\widehat{G}(p)\widehat{g}(p)}{p^{2}-ae^{-iph}}, \quad
p^{2}\widehat{u_{m}}(p)=\sqrt{2\pi}\frac{p^{2}\widehat{G_{m}}(p)\widehat{g_{m}}(p)}{p^{2}-ae^{-iph}}, \quad m\in {\mathbb N}.
\end{equation}
Obviously,
$$
|p^{2}\widehat{u_{m}}(p)-p^{2}\widehat{u}(p)|\leq \sqrt{2\pi}\Bigg\|\frac{p^{2}\widehat{G_{m}}(p)}{p^{2}-ae^{-iph}}-
\frac{p^{2}\widehat{G}(p)}{p^{2}-ae^{-iph}}\Bigg\|_{L^{\infty}({\mathbb R})}|\widehat{g}(p)|+
$$
$$
\sqrt{2\pi}\Bigg\|\frac{p^{2}\widehat{G_{m}}(p)}{p^{2}-ae^{-iph}}\Bigg\|_{L^{\infty}({\mathbb R})}|\widehat{g_{m}}(p)-\widehat{g}(p)|.
$$
Using bound (\ref{fumuxn}), we obtain
$$
\Bigg\|\frac{d^{2}u_{m}}{dx^{2}}-\frac{d^{2}u}{dx^{2}}\Bigg\|_{L^{2}({\mathbb R})}\leq \sqrt{2\pi}\Bigg\|\frac{p^{2}\widehat{G_{m}}(p)}{p^{2}-
ae^{-iph}}-\frac{p^{2}\widehat{G}(p)}{p^{2}-ae^{-iph}}\Bigg\|_{L^{\infty}({\mathbb R})}\|F(u(x), x)\|_{L^{2}({\mathbb R})}+
$$
$$
\sqrt{2\pi}\Bigg\|\frac{p^{2}\widehat{G_{m}}(p)}{p^{2}-ae^{-iph}}\Bigg\|_{L^{\infty}({\mathbb R})}l
\|u_{m}(x)-u(x)\|_{L^{2}({\mathbb R})}.
$$
Let us recall the result of Lemma A2 along with statement (\ref{umtoul2}). This means that
\begin{equation}
\label{umtoul22}
\frac{d^{2}u_{m}}{dx^{2}}\to \frac{d^{2}u}{dx^{2}}, \quad m\to \infty
\end{equation}
in $L^{2}({\mathbb R})$. By means of definition (\ref{sob}) of the norm,  propositions (\ref{umtoul2}) and (\ref{umtoul22}), we have
$$
u_{m}(x)\to u(x), \quad m\to \infty
$$
in $H^{2}({\mathbb R})$.

Suppose that for some $m\in {\mathbb N}$ the unique solution of (\ref{id1m}) considered above vanishes identically on the real line. This will contradict to the given condition that the Fourier images of $G_{m}(x)$ and $F(0, x)$ are nontrivial on a
set of nonzero Lebesgue measure in ${\mathbb R}$.
The analogous argument is valid for the unique solution $u(x)$ of limiting
problem  (\ref{id1}).   \hfill\lanbox


\setcounter{equation}{0}

\section{Appendix}

\bigskip

Consider a function  $G(x): {\mathbb R}\to {\mathbb R}$.
Let us denote its standard Fourier transform using the hat symbol as
\begin{equation}
\label{f}
{\cal F}[G(x)](p)=\widehat{G}(p):=\frac{1}{\sqrt{2\pi}}\int_{-\infty}^{\infty}G(x)
e^{-ipx}dx, \ p\in {\mathbb R}.
\end{equation}
Clearly,
\begin{equation}
\label{inf1}
\|\widehat{G}(p)\|_{L^{\infty}({\mathbb R})}\leq \frac{1}{\sqrt{2\pi}}
\|G(x)\|_{L^{1}({\mathbb R})}
\end{equation}
and
$\displaystyle{G(x)=\frac{1}{\sqrt{2\pi}}\int_{-\infty}^{\infty}
\widehat{G}(q)e^{iqx}dq, \ x\in {\mathbb R}.}$
Let us introduce the technical expressions to address the solvability of equations (\ref{id1m}) with $m\in {\mathbb N}$ as
\begin{equation}
\label{Na}
N_{a, \ h, \ m}:=max\Big\{
\Big\|\frac{\widehat{G}_{m}(p)}{p^{2}-ae^{-iph}}\Big\|_{L^{\infty}({\mathbb R})},
\quad
\Big\|\frac{p^{2}\widehat{G}_{m}(p)}{p^{2}-ae^{-iph}}\Big\|_{L^{\infty}({\mathbb R})}
\Big\}.
\end{equation}
Analogously, in the limiting case
\begin{equation}
\label{Nal}
N_{a, \ h}:=max\Big\{
\Big\|\frac{\widehat{G}(p)}{p^{2}-ae^{-iph}}\Big\|_{L^{\infty}({\mathbb R})},
\quad
\Big\|\frac{p^{2}\widehat{G}(p)}{p^{2}-ae^{-iph}}\Big\|_{L^{\infty}({\mathbb R})}
\Big\}.
\end{equation}
In formulas (\ref{Na}) and (\ref{Nal}) above we have the constants
$a>0, \ h\in {\mathbb R}, \ h\neq 0$.

\bigskip

\noindent
{\bf Lemma A1.} {\it Let $G(x): {\mathbb R}\to {\mathbb R}$ and
$G(x)\in L^{1}({\mathbb R})$.

\noindent
a) If $\displaystyle{h\neq \frac{2\pi n}{\sqrt{a}}, \ n\in {\mathbb Z}}$ then
$N_{a, \ h}<\infty$.

\noindent
b) When
$\displaystyle{h=\frac{2\pi n}{\sqrt{a}}, \ n\in {\mathbb Z}, \ n\neq 0}$,
let in addition $xG(x)\in L^{1}({\mathbb R})$. Then $N_{a, \ h}<\infty$ if
and only if the orthogonality conditions
\begin{equation}
\label{oc}
\Bigg(G(x), \frac{e^{\pm i\sqrt{a}x}}{\sqrt{2\pi}} \Bigg)_{L^{2}({\mathbb R})}=0
\end{equation}
hold.}

\bigskip

\noindent
{\it Proof.} First of all, it can be easily verified that
$\displaystyle{\frac{\widehat{G}(p)}{p^{2}-ae^{-iph}}}$ is bounded if and only
if
$\displaystyle{\frac{p^{2}\widehat{G}(p)}{p^{2}-ae^{-iph}}}$ belongs to
$L^{\infty}({\mathbb R})$.
Indeed, we have the identity
\begin{equation}
\label{p2Ghpeq}
\frac{p^{2}\widehat{G}(p)}{p^{2}-ae^{-iph}}=\widehat{G}(p)+
\frac{ae^{-iph}\widehat{G}(p)}{p^{2}-ae^{-iph}}.
\end{equation}
Let us suppose that
$\displaystyle{\frac{\widehat{G}(p)}{p^{2}-ae^{-iph}}\in L^{\infty}({\mathbb R})}$.
Using (\ref{p2Ghpeq}) along with (\ref{inf1}) and the given assumptions, we
arrive at
$$
\Bigg\|\frac{p^{2}\widehat{G}(p)}{p^{2}-ae^{-iph}}\Bigg\|_{L^{\infty}({\mathbb R})}\leq
\|\widehat{G}(p)\|_{L^{\infty}({\mathbb R})}+
a\Bigg\|\frac{\widehat{G}(p)}{p^{2}-ae^{-iph}}\Bigg\|_{L^{\infty}({\mathbb R})}<\infty,
$$
such that
\begin{equation}
\label{p2Ghp}
\frac{p^{2}\widehat{G}(p)}{p^{2}-ae^{-iph}}\in L^{\infty}({\mathbb R})
\end{equation}
as well. On the other hand, if (\ref{p2Ghp}) holds, then by means of
equality (\ref{p2Ghpeq}) and upper bound (\ref{inf1}) under our assumptions
$$
\Bigg\|\frac{\widehat{G}(p)}{p^{2}-ae^{-iph}}\Bigg\|_{L^{\infty}({\mathbb R})}\leq
\frac{1}{a}
\Bigg\|\frac{p^{2}\widehat{G}(p)}{p^{2}-ae^{-iph}}\Bigg\|_{L^{\infty}({\mathbb R})}+
\frac{1}{a}\|\widehat{G}(p)\|_{L^{\infty}({\mathbb R})}<\infty.
$$
Hence,
$\displaystyle{\frac{\widehat{G}(p)}{p^{2}-ae^{-iph}}}$ belongs to
$L^{\infty}({\mathbb R})$.

Let us first consider the case when $\displaystyle{h\neq \frac{2\pi n}{\sqrt{a}}, \ n\in {\mathbb Z}}$. Recall inequalities (\ref{alfa})
and (\ref{inf1}). Hence,
$$
\Bigg|\frac{\widehat{G}(p)}{p^{2}-ae^{-iph}}\Bigg|\leq \frac{|\widehat{G}(p)|}{\sqrt{\alpha}}\leq
\frac{1}{\sqrt{\alpha}}\frac{1}{\sqrt{2\pi}}\|G(x)\|_{L^{1}({\mathbb R})},
$$
such that
$$
\frac{\widehat{G}(p)}{p^{2}-ae^{-iph}}\in L^{\infty}({\mathbb R}).
$$
According to the reasoning above, $N_{a, \ h}<\infty$ in the situation a) of our lemma.

Then we cover the case when $\displaystyle{h=\frac{2\pi n}{\sqrt{a}}, \ n\in {\mathbb Z}, \ n\neq 0}$. Let us write
$$
\frac{{\widehat G}(p)}{p^{2}-ae^{-iph}}=\frac{{\widehat G}(p)}{p^{2}-ae^{-iph}}\chi_{I_{\delta}^{+}}+
\frac{{\widehat G}(p)}{p^{2}-ae^{-iph}}\chi_{I_{\delta}^{-}}+
$$
\begin{equation}
\label{uhpf4g}
\frac{{\widehat G}(p)}{p^{2}-ae^{-iph}}\chi_{I_{\delta}^{c +}}+
\frac{{\widehat G}(p)}{p^{2}-ae^{-iph}}\chi_{I_{\delta}^{c -}}.
\end{equation}
The third term in the right side of (\ref{uhpf4g}) can be estimated from above in the absolute value using (\ref{inf1}) as
$$
\frac{|{\widehat G}(p)|\chi_{I_{\delta}^{c +}}}{\sqrt{(p^{2}-a)^{2}+2ap^{2}(1-cos(ph))}}\leq
\frac{|{\widehat G}(p)|\chi_{I_{\delta}^{c +}}}{\sqrt{(p-\sqrt{a})^{2}(p+\sqrt{a})^{2}}}\leq
$$
$$
\frac{|{\widehat G}(p)|}{\delta \sqrt{a}}\leq \frac{1}{\sqrt{2\pi a}\delta}\|G(x)\|_{L^{1}({\mathbb R})}<\infty
$$
as assumed, so that
$$
\frac{{\widehat G}(p)}{p^{2}-ae^{-iph}}\chi_{I_{\delta}^{c +}}\in L^{\infty}({\mathbb R}).
$$
The fourth term in the right side of (\ref{uhpf4g}) can be bounded from above in the absolute value via (\ref{inf1}) as
$$
\frac{|{\widehat G}(p)|\chi_{I_{\delta}^{c -}}}{\sqrt{(p^{2}-a)^{2}+2ap^{2}(1-cos(ph))}}\leq
\frac{|{\widehat G}(p)|\chi_{I_{\delta}^{c -}}}{\sqrt{(p-\sqrt{a})^{2}(p+\sqrt{a})^{2}}}\leq
$$
$$
\frac{|{\widehat G}(p)|}{\delta \sqrt{a}}\leq \frac{1}{\sqrt{2\pi a}\delta}\|G(x)\|_{L^{1}({\mathbb R})}<\infty
$$
as above, such that
$$
\frac{{\widehat G}(p)}{p^{2}-ae^{-iph}}\chi_{I_{\delta}^{c -}}\in L^{\infty}({\mathbb R}).
$$
Evidently,
$$
{\widehat G}(p)={\widehat G}(\sqrt{a})+\int_{\sqrt{a}}^{p}\frac{d{\widehat G}(s)}{ds}ds.
$$
Thus, the first term in the right side of (\ref{uhpf4g}) is equal to
\begin{equation}
\label{ghp2aidp}
\frac{{\widehat G}(\sqrt{a})}{p^{2}-ae^{-iph}}\chi_{I_{\delta}^{+}}+
\frac{\int_{\sqrt{a}}^{p}\frac{d{\widehat G}(s)}{ds}ds}{p^{2}-ae^{-iph}}\chi_{I_{\delta}^{+}}.
\end{equation}
By means of the analog of inequality (\ref{dfhpdp}), we arrive at
$$
\Bigg|\frac{\int_{\sqrt{a}}^{p}\frac{d{\widehat G}(s)}{ds}ds}{p^{2}-ae^{-iph}}\chi_{I_{\delta}^{+}}\Bigg|\leq \frac
{\|xG(x)\|_{L^{1}({\mathbb R})}|p-\sqrt{a}|\chi_{I_{\delta}^{+}}}{\sqrt{2\pi}\sqrt{(p^{2}-a)^{2}+2ap^{2}(1-\cos(ph))}}\leq
$$
$$
\frac{1}{\sqrt{2\pi}}\frac{\|xG(x)\|_{L^{1}({\mathbb R})}\chi_{I_{\delta}^{+}}}{\sqrt{(p+\sqrt{a})^{2}}}\leq
\frac{1}{\sqrt{2\pi a}}\|xG(x)\|_{L^{1}({\mathbb R})}<\infty
$$
as we assume. Hence,
$$
\frac{\int_{\sqrt{a}}^{p}\frac{d{\widehat G}(s)}{ds}ds}{p^{2}-ae^{-iph}}\chi_{I_{\delta}^{+}}\in L^{\infty}({\mathbb R}).
$$
Consider the remaining term in (\ref{ghp2aidp}) given by
\begin{equation}
\label{ghraidp}
\frac{{\widehat G}(\sqrt{a})}{p^{2}-ae^{-iph}}\chi_{I_{\delta}^{+}}.
\end{equation}
Its absolute value equals to
$$
\frac{|{\widehat G}(\sqrt{a})|\chi_{I_{\delta}^{+}}}{\sqrt{(p^{2}-a)^{2}+2ap^{2}\Big(1-\cos\Big(p\frac{2\pi n}{\sqrt{a}}\Big)\Big)}}.
$$
Thus, expression (\ref{ghraidp}) is bounded on the real line if and only if ${\widehat G}(\sqrt{a})=0$. This is equivalent to the
orthogonality condition
\begin{equation}
\label{ocg1}
\Bigg(G(x), \frac{e^{i\sqrt{a}x}}{\sqrt{2\pi}} \Bigg)_{L^{2}({\mathbb R})}=0.
\end{equation}
Clearly,
$$
{\widehat G}(p)={\widehat G}(-\sqrt{a})+\int_{-\sqrt{a}}^{p}\frac{d{\widehat G}(s)}{ds}ds.
$$
Hence, the second term in the right side of (\ref{uhpf4g}) is equal to
\begin{equation}
\label{ghp2aidp2}
\frac{{\widehat G}(-\sqrt{a})}{p^{2}-ae^{-iph}}\chi_{I_{\delta}^{-}}+
\frac{\int_{-\sqrt{a}}^{p}\frac{d{\widehat G}(s)}{ds}ds}{p^{2}-ae^{-iph}}\chi_{I_{\delta}^{-}}.
\end{equation}
By virtue of the analog of bound (\ref{dfhpdp}), we derive
$$
\Bigg|\frac{\int_{-\sqrt{a}}^{p}\frac{d{\widehat G}(s)}{ds}ds}{p^{2}-ae^{-iph}}\chi_{I_{\delta}^{-}}\Bigg|\leq \frac
{\|xG(x)\|_{L^{1}({\mathbb R})}|p+\sqrt{a}|\chi_{I_{\delta}^{-}}}{\sqrt{2\pi}\sqrt{(p^{2}-a)^{2}+2ap^{2}(1-\cos(ph))}}\leq
$$
$$
\frac{1}{\sqrt{2\pi}}\frac{\|xG(x)\|_{L^{1}({\mathbb R})}\chi_{I_{\delta}^{-}}}{\sqrt{(p-\sqrt{a})^{2}}}\leq
\frac{1}{\sqrt{2\pi a}}\|xG(x)\|_{L^{1}({\mathbb R})}<\infty
$$
as above. Thus,
$$
\frac{\int_{-\sqrt{a}}^{p}\frac{d{\widehat G}(s)}{ds}ds}{p^{2}-ae^{-iph}}\chi_{I_{\delta}^{-}}\in L^{\infty}({\mathbb R}).
$$
The remaining term in (\ref{ghp2aidp2}) is
\begin{equation}
\label{ghraidp}
\frac{{\widehat G}(-\sqrt{a})}{p^{2}-ae^{-iph}}\chi_{I_{\delta}^{-}}.
\end{equation}
Its absolute value is given by
$$
\frac{|{\widehat G}(-\sqrt{a})|\chi_{I_{\delta}^{-}}}{\sqrt{(p^{2}-a)^{2}+2ap^{2}\Big(1-\cos\Big(p\frac{2\pi n}{\sqrt{a}}\Big)\Big)}}.
$$
Obviously, expression (\ref{ghraidp}) belongs to $L^{\infty}({\mathbb R})$ if and only if ${\widehat G}(-\sqrt{a})=0$.
This is equivalent to the orthogonality relation
\begin{equation}
\label{ocg2}
\Bigg(G(x), \frac{e^{-i\sqrt{a}x}}{\sqrt{2\pi}} \Bigg)_{L^{2}({\mathbb R})}=0.
\end{equation}
Therefore, by means of formulas (\ref{ocg1}) and (\ref{ocg2}) we obtain that $N_{a, \ h}<\infty$ if
and only if orthogonality conditions (\ref{oc}) are valid in the situation when
$\displaystyle{h=\frac{2\pi n}{\sqrt{a}}, \ n\in {\mathbb Z}, \ n\neq 0}$.  \hfill\lanbox

\bigskip

Note that in the case a) of the lemma above the orthogonality relations are
not needed. But in the case b) the orthogonality conditions are analogous
to the ones we have in the case a) of Lemma A1 of ~\cite{VV1021}.
The final auxiliary statement of the article is as follows.

\bigskip

\noindent
{\bf Lemma A2.} {\it Let $m\in {\mathbb N}$, the kernels
$G_{m}(x): {\mathbb R}\to {\mathbb R}$ and $G_{m}(x)\in L^{1}({\mathbb R})$,
such that $G_{m}(x)\to G(x)$ in $L^{1}({\mathbb R})$ as $m\to \infty$.

\noindent
a) If $\displaystyle{h\neq \frac{2\pi n}{\sqrt{a}}, \ n\in {\mathbb Z}}$,  we assume that
\begin{equation}
\label{nahm}
2\sqrt{\pi}N_{a, \ h, \ m}l\leq 1-\varepsilon, \quad a>0, \quad h\in {\mathbb R},
\quad h\neq 0
\end{equation}
holds for all $m\in {\mathbb N}$ for a certain fixed $0<\varepsilon<1$.

\noindent
b) When $\displaystyle{h=\frac{2\pi n}{\sqrt{a}}, \ n\in {\mathbb Z}, \ n\neq 0}$,
let $xG_{m}(x)\in L^{1}({\mathbb R})$, so that
$xG_{m}(x)\to xG(x)$ in $L^{1}({\mathbb R})$ as $m\to \infty$ and in addition
\begin{equation}
\label{ocm}
\Bigg(G_{m}(x), \frac{e^{\pm i\sqrt{a}x}}{\sqrt{2\pi}} \Bigg)_{L^{2}({\mathbb R})}=0,
\quad m\in {\mathbb N}.
\end{equation}
Moreover, we assume that inequality (\ref{nahm}) is valid
for all $m\in {\mathbb N}$ for some fixed $0<\varepsilon<1$. Then
\begin{equation}
\label{conv1}
\frac{\widehat{G}_{m}(p)}{p^{2}-ae^{-iph}}\to
\frac{\widehat{G}(p)}{p^{2}-ae^{-iph}}, \quad
\frac{p^{2}\widehat{G}_{m}(p)}{p^{2}-ae^{-iph}}\to
\frac{p^{2}\widehat{G}(p)}{p^{2}-ae^{-iph}}
\end{equation}
in $L^{\infty}({\mathbb R})$ as $m\to \infty$, so that
\begin{equation}
\label{conv2}
\Bigg\|\frac{\widehat{G}_{m}(p)}{p^{2}-ae^{-iph}}\Bigg\|_{L^{\infty}({\mathbb R})}\to
\Bigg\|\frac{\widehat{G}(p)}{p^{2}-ae^{-iph}}\Bigg\|_{L^{\infty}({\mathbb R})}, \quad
m\to \infty,
\end{equation}
\begin{equation}
\label{conv3}
\Bigg\|\frac{p^{2}\widehat{G}_{m}(p)}{p^{2}-ae^{-iph}}\Bigg\|_{L^{\infty}({\mathbb R})}\to
\Bigg\|\frac{p^{2}\widehat{G}(p)}{p^{2}-ae^{-iph}}\Bigg\|_{L^{\infty}({\mathbb R})}, \quad
m\to \infty.
\end{equation}
Furthermore,
\begin{equation}
\label{nah}
2\sqrt{\pi}N_{a, \ h}l\leq 1-\varepsilon, \quad a>0, \quad h\in {\mathbb R},
\quad h\neq 0
\end{equation}
is valid.}

\bigskip

\noindent
{\it Proof.} Recall inequality (\ref{inf1}). Hence,
\begin{equation}
\label{gmhgph}
\|\widehat{G}_{m}(p)-\widehat{G}(p)\|_{L^{\infty}({\mathbb R})}\leq \frac{1}{\sqrt{2\pi}}\|G_{m}(x)-G(x)\|_{L^{1}({\mathbb R})}\to 0, \quad
m\to \infty
\end{equation}
as we assume. Clearly, (\ref{conv2}) and (\ref{conv3}) follow from (\ref{conv1}) by means of the
standard triangle inequality.  It can easily verified that in (\ref{conv1}) the first statement yields the second one. Indeed,  let us express
$$
\frac{p^{2}\widehat{G}_{m}(p)}{p^{2}-ae^{-iph}}-\frac{p^{2}\widehat{G}(p)}{p^{2}-ae^{-iph}}=\widehat{G}_{m}(p)-\widehat{G}(p)+
ae^{-iph}\Bigg[\frac{\widehat{G}_{m}(p)}{p^{2}-ae^{-iph}}-\frac{\widehat{G}(p)}{p^{2}-ae^{-iph}}\Bigg],
$$
so that
$$
\Bigg\|\frac{p^{2}\widehat{G}_{m}(p)}{p^{2}-ae^{-iph}}-\frac{p^{2}\widehat{G}(p)}{p^{2}-ae^{-iph}}\Bigg\|_{L^{\infty}({\mathbb R})}\leq
\|\widehat{G}_{m}(p)-\widehat{G}(p)\|_{L^{\infty}({\mathbb R})}+
$$
\begin{equation}
\label{p2gmhpghpn}
a\Bigg\|\frac{\widehat{G}_{m}(p)}{p^{2}-ae^{-iph}}-\frac{\widehat{G}(p)}{p^{2}-ae^{-iph}}\Bigg\|_{L^{\infty}({\mathbb R})}.
\end{equation}
The second term in the right side of  (\ref{p2gmhpghpn}) converges to zero as $m\to \infty$ according to the first proposition in
 (\ref{conv1}). The first term in the right side of  (\ref{p2gmhpghpn}) tends to zero as $m\to \infty$ by virtue of  (\ref{gmhgph}).
This means that
$$
\frac{p^{2}\widehat{G}_{m}(p)}{p^{2}-ae^{-iph}}\to \frac{p^{2}\widehat{G}(p)}{p^{2}-ae^{-iph}}, \quad m\to \infty
$$
in $L^{\infty}({\mathbb R})$.

First we address the case when $\displaystyle{h\neq \frac{2\pi n}{\sqrt{a}}, \ n\in {\mathbb Z}}$ using bounds (\ref{alfa}) and (\ref{inf1}).
Thus,
$$
\Bigg|\frac{\widehat{G}_{m}(p)}{p^{2}-ae^{-iph}}-\frac{\widehat{G}(p)}{p^{2}-ae^{-iph}}\Bigg|\leq \frac{\|G_{m}(x)-G(x)\|_{L^{1}({\mathbb R})}}
{\sqrt{2\pi}\sqrt{(p^{2}-a)^{2}+2ap^{2}(1-\cos(ph))}}\leq
$$
$$
\frac{1}{\sqrt{2\pi \alpha}}\|G_{m}(x)-G(x)\|_{L^{1}({\mathbb R})},
$$
such that
$$
\Bigg\|\frac{\widehat{G}_{m}(p)}{p^{2}-ae^{-iph}}-\frac{\widehat{G}(p)}{p^{2}-ae^{-iph}}\Bigg\|_{L^{\infty}({\mathbb R})} \leq
\frac{1}{\sqrt{2\pi \alpha}}\|G_{m}(x)-G(x)\|_{L^{1}({\mathbb R})}\to 0
$$
as $m\to \infty$ due to the given condition. Therefore,
$$
\frac{\widehat{G}_{m}(p)}{p^{2}-ae^{-iph}}\to \frac{\widehat{G}(p)}{p^{2}-ae^{-iph}}, \quad m\to \infty
$$
in $L^{\infty}({\mathbb R})$ in the situation a) of our lemma.

Finally, we discuss the case when
$\displaystyle{h=\frac{2\pi n}{\sqrt{a}}, \ n\in {\mathbb Z}, \ n\neq 0}$. Obviously,
$$
\frac{\widehat{G}_{m}(p)}{p^{2}-ae^{-iph}}-\frac{\widehat{G}(p)}{p^{2}-ae^{-iph}}=\frac{\widehat{G}_{m}(p)-\widehat{G}(p)}{p^{2}-ae^{-iph}}\chi_{I_{\delta}^{+}}+
\frac{\widehat{G}_{m}(p)-\widehat{G}(p)}{p^{2}-ae^{-iph}}\chi_{I_{\delta}^{-}}+
$$
\begin{equation}
\label{uhpf4gd}
\frac{\widehat{G}_{m}(p)-\widehat{G}(p)}{p^{2}-ae^{-iph}}\chi_{I_{\delta}^{c +}}+
\frac{\widehat{G}_{m}(p)-\widehat{G}(p)}{p^{2}-ae^{-iph}}\chi_{I_{\delta}^{c -}}.
\end{equation}
Let us use inequality (\ref{inf1}) to estimate the third term in the right side of  (\ref{uhpf4gd}) from above in the absolute value as
$$
\frac{|\widehat{G}_{m}(p)-\widehat{G}(p)|\chi_{I_{\delta}^{c +}}}{\sqrt{(p^{2}-a)^{2}+2ap^{2}(1-\cos(ph))}}\leq
\frac{|\widehat{G}_{m}(p)-\widehat{G}(p)|\chi_{I_{\delta}^{c +}}}{|p-\sqrt{a}||p+\sqrt{a}|}\leq
$$
$$
\frac{|\widehat{G}_{m}(p)-\widehat{G}(p)|}{\delta \sqrt{a}}\leq \frac{\|G_{m}(x)-G(x)\|_{L^{1}({\mathbb R})}}{\sqrt{2\pi a}\delta},
$$
so that
$$
\Bigg\|\frac{\widehat{G}_{m}(p)-\widehat{G}(p)}{p^{2}-ae^{-iph}}\chi_{I_{\delta}^{c +}}\Bigg\|_{L^{\infty}({\mathbb R})}\leq
\frac{\|G_{m}(x)-G(x)\|_{L^{1}({\mathbb R})}}{\sqrt{2\pi a}\delta}\to 0, \quad m\to \infty
$$
as assumed. Similarly, the fourth term in the right side of  (\ref{uhpf4gd}) can be bounded from above in the absolute value via
 (\ref{inf1})  as
$$
\frac{|\widehat{G}_{m}(p)-\widehat{G}(p)|\chi_{I_{\delta}^{c -}}}{\sqrt{(p^{2}-a)^{2}+2ap^{2}(1-cos(ph))}}\leq
\frac{|\widehat{G}_{m}(p)-\widehat{G}(p)|\chi_{I_{\delta}^{c -}}}{|p-\sqrt{a}||p+\sqrt{a}|}\leq
$$
$$
\frac{|\widehat{G}_{m}(p)-\widehat{G}(p)|}{\delta \sqrt{a}}\leq \frac{\|G_{m}(x)-G(x)\|_{L^{1}({\mathbb R})}}{\sqrt{2\pi a}\delta},
$$
such that
$$
\Bigg\|\frac{\widehat{G}_{m}(p)-\widehat{G}(p)}{p^{2}-ae^{-iph}}\chi_{I_{\delta}^{c -}}\Bigg\|_{L^{\infty}({\mathbb R})}\leq
\frac{\|G_{m}(x)-G(x)\|_{L^{1}({\mathbb R})}}{\sqrt{2\pi a}\delta}\to 0, \quad m\to \infty
$$
as above. Let us recall orthogonality conditions (\ref{ocm}). It can be trivially checked that the analogous relations will hold in the limit.
Obviously,
$$
\Bigg|\Bigg(G(x), \frac{e^{\pm i\sqrt{a}x}}{\sqrt{2\pi}}\Bigg)_{L^{2}({\mathbb R})}\Bigg|=
\Bigg|\Bigg(G(x)-G_{m}(x), \frac{e^{\pm i\sqrt{a}x}}{\sqrt{2\pi}}\Bigg)_{L^{2}({\mathbb R})}\Bigg|\leq
$$
$$
\frac{1}{\sqrt{2\pi}}\|G_{m}(x)-G(x)\|_{L^{1}({\mathbb R})}\to 0, \quad m\to \infty
$$
via the given condition. This means that
\begin{equation}
\label{ocml}
\Bigg(G(x), \frac{e^{\pm i\sqrt{a}x}}{\sqrt{2\pi}}\Bigg)_{L^{2}({\mathbb R})}=0.
\end{equation}
By virtue of  (\ref{ocml}) and  (\ref{ocm}), we have
$$
\widehat{G}(\sqrt{a})=0, \quad \widehat{G}_{m}(\sqrt{a})=0, \quad m\in {\mathbb N}.
$$
Therefore,
$$
\widehat{G}(p)=\int_{\sqrt{a}}^{p}\frac{d\widehat{G}(s)}{ds}ds, \quad
\widehat{G}_{m}(p)=\int_{\sqrt{a}}^{p}\frac{d\widehat{G}_{m}(s)}{ds}ds, \quad m\in {\mathbb N}.
$$
Then  the first term in the right side of  (\ref{uhpf4gd}) is given by
\begin{equation}
\label{rapgmghs}
\frac{\int_{\sqrt{a}}^{p}\Big[\frac{d\widehat{G}_{m}(s)}{ds}-\frac{d\widehat{G}(s)}{ds}\Big]ds}{p^{2}-ae^{-iph}}\chi_{I_{\delta}^{+}}.
\end{equation}
Let us use the analog of inequality (\ref{dfmhfhm}) to derive the upper bound on (\ref{rapgmghs}) in the absolute value
as
$$
\Bigg|\frac{\int_{\sqrt{a}}^{p}\Big[\frac{d\widehat{G}_{m}(s)}{ds}-\frac{d\widehat{G}(s)}{ds}\Big]ds}{p^{2}-ae^{-iph}}\chi_{I_{\delta}^{+}}\Bigg|\leq \frac{\|xG_{m}(x)-xG(x)\|_{L^{1}({\mathbb R})}|p-\sqrt{a}|}{\sqrt{2\pi}
\sqrt{(p^{2}-a)^{2}+2ap^{2}(1-\cos(ph))}}\chi_{I_{\delta}^{+}}\leq
$$
$$
\frac{\|xG_{m}(x)-xG(x)\|_{L^{1}({\mathbb R})}}{\sqrt{2\pi}\sqrt{(p+\sqrt{a})^{2}}}\chi_{I_{\delta}^{+}}\leq
\frac{\|xG_{m}(x)-xG(x)\|_{L^{1}({\mathbb R})}}{\sqrt{2\pi a}}\chi_{I_{\delta}^{+}},
$$
so that
$$
\Bigg\|\frac{\int_{\sqrt{a}}^{p}\Big[\frac{d\widehat{G}_{m}(s)}{ds}-\frac{d\widehat{G}(s)}{ds}\Big]ds}{p^{2}-ae^{-iph}}\chi_{I_{\delta}^{+}}
\Bigg\|_{L^{\infty}({\mathbb R})}\leq \frac{\|xG_{m}(x)-xG(x)\|_{L^{1}({\mathbb R})}}{\sqrt{2\pi a}}\to 0
$$
as $m\to \infty$ as we assume.

Recall (\ref{ocml}) and  (\ref{ocm}). Hence,
$$
\widehat{G}(-\sqrt{a})=0, \quad \widehat{G}_{m}(-\sqrt{a})=0, \quad m\in {\mathbb N}.
$$
This means that
$$
\widehat{G}(p)=\int_{-\sqrt{a}}^{p}\frac{d\widehat{G}(s)}{ds}ds, \quad
\widehat{G}_{m}(p)=\int_{-\sqrt{a}}^{p}\frac{d\widehat{G}_{m}(s)}{ds}ds, \quad m\in {\mathbb N}
$$
and the second term in the right side of  (\ref{uhpf4gd}) equals to
\begin{equation}
\label{rapgmghsm}
\frac{\int_{-\sqrt{a}}^{p}\Big[\frac{d\widehat{G}_{m}(s)}{ds}-\frac{d\widehat{G}(s)}{ds}\Big]ds}{p^{2}-ae^{-iph}}\chi_{I_{\delta}^{-}}.
\end{equation}
By virtue of the analog of (\ref{dfmhfhm}) we obtain the estimate from above on  expression (\ref{rapgmghsm}) in the absolute value
as
$$
\Bigg|\frac{\int_{-\sqrt{a}}^{p}\Big[\frac{d\widehat{G}_{m}(s)}{ds}-\frac{d\widehat{G}(s)}{ds}\Big]ds}{p^{2}-ae^{-iph}}\chi_{I_{\delta}^{-}}\Bigg|\leq \frac{\|xG_{m}(x)-xG(x)\|_{L^{1}({\mathbb R})}|p+\sqrt{a}|}{\sqrt{2\pi}
\sqrt{(p^{2}-a)^{2}+2ap^{2}(1-\cos(ph))}}\chi_{I_{\delta}^{-}}\leq
$$
$$
\frac{\|xG_{m}(x)-xG(x)\|_{L^{1}({\mathbb R})}}{\sqrt{2\pi}\sqrt{(p-\sqrt{a})^{2}}}\chi_{I_{\delta}^{-}}\leq
\frac{\|xG_{m}(x)-xG(x)\|_{L^{1}({\mathbb R})}}{\sqrt{2\pi a}}\chi_{I_{\delta}^{-}},
$$
such that
$$
\Bigg\|\frac{\int_{-\sqrt{a}}^{p}\Big[\frac{d\widehat{G}_{m}(s)}{ds}-\frac{d\widehat{G}(s)}{ds}\Big]ds}{p^{2}-ae^{-iph}}\chi_{I_{\delta}^{-}}
\Bigg\|_{L^{\infty}({\mathbb R})}\leq \frac{\|xG_{m}(x)-xG(x)\|_{L^{1}({\mathbb R})}}{\sqrt{2\pi a}}\to 0
$$
as $m\to \infty$ due to the given condition.  Thus,
$$
\frac{\widehat{G}_{m}(p)}{p^{2}-ae^{-iph}}\to \frac{\widehat{G}(p)}{p^{2}-ae^{-iph}}, \quad m\to \infty
$$
in $L^{\infty}({\mathbb R})$ in the case b) of the lemma as well. Obviously, under the stated assumptions
$$
N_{a, \ h}<\infty, \quad N_{a, \ h, \ m}<\infty, \quad m\in {\mathbb N}, \quad a>0, \quad h\in {\mathbb R}, \quad h\neq 0
$$
according to the result of Lemma A1. Recall  (\ref{nahm}).  A trivial limiting argument relying on (\ref{conv2}) and (\ref{conv3}) yields
(\ref{nah}).  \hfill\lanbox


\section*{Acknowledgements}

The first author is grateful to Israel Michael Sigal for the partial support by the NSERC
grant NA 7901.
The second author has been supported by the RUDN
University Strategic Academic Leadership Program.


\bigskip

\begin{thebibliography}{ccccccccc}

\medskip


\bibitem{Ag}
M.S. Agranovich. {\em Elliptic boundary problems.}
Encyclopaedia Math. Sci., vol. 79, Partial Differential Equations,
IX, Springer, Berlin (1997), pp. 1--144

\medskip

\bibitem{AMP14} G.L. ~Alfimov, E.V. ~Medvedeva, D.E. ~Pelinovsky.
{\em Wave systems with an infinite number of localized traveling waves,}
Phys. Rev. Lett., {\bf 112} (2014), 054103, 5 pp.


\medskip

\bibitem{AKLP19}
G.L. ~Alfimov, A.S. ~Korobeinikov, C.J. ~Lustri, D.E. ~Pelinovsky.
{\em Standing lattice solitons in the discrete NLS equation with saturation,}
Nonlinearity, {\bf  32} (2019), no. 9, 3445--3484

\medskip

\bibitem{ABVV10}
 N. ~Apreutesei, N. ~Bessonov, V. ~Volpert, V. ~Vougalter.
 {\em Spatial structures and generalized travelling waves for an integro-
 differential equation,}  Discrete Contin. Dyn. Syst., Ser. B, {\bf 13}
(2010), no. 3, 537--557

\medskip

\bibitem{BNPR09}
H. ~Berestycki, G. ~Nadin, B. ~Perthame, L. ~Ryzhik.
{\em The non-local Fisher-KPP equation: traveling waves and steady states,}
Nonlinearity, {\bf 22} (2009), no. 12, 2813--2844

\medskip

\bibitem{BO86}
H. ~Brezis, L. ~Oswald. {\em Remarks on sublinear elliptic equations,}
Nonlinear Anal., Theory Methods Appl., {\bf 10} (1986), 55--64

\medskip

\bibitem{DMV05} A.~Ducrot, M. ~Marion,  V. ~Volpert. {\em Systemes de
r\'eaction-diffusion sans propri\'et\'e de Fredholm,} C. R. Math. Acad. Sci. Paris,  {\bf
340} (2005), no. 9, 659--664

\medskip

\bibitem{DMV08}
 A. ~Ducrot, M. ~Marion,  V. ~Volpert.
 {\em Reaction-diffusion problems with non-Fredholm operators,}
 Adv. Differ. Equ., {\bf 13} (2008), no. 11--12, 1151--1192

\medskip

\bibitem{DMV08b}
 A. ~Ducrot, M. ~Marion, V. ~Volpert.
 {\em Reaction-diffusion waves
(with the Lewis number different from 1).} Mathematics and Mathematical Modelling. Paris: Publibooks (2009), 113 pp.

\medskip

\bibitem{E09}
M. ~Efendiev. {\em Fredholm structures, topological invariants and
applications,} AIMS Series on Differential Equations \& Dynamical Systems,
{\bf 3}. American Institute of Mathematical Sciences (AIMS), Springfield, MO
(2009), 205 pp.

\medskip

\bibitem{E18}
M. ~Efendiev. {\em Symmetrization and stabilization of solutions of nonlinear elliptic equations,}
Fields Institute Monographs,  {\bf 36}.  The Fields Institute for Research in the Mathematical Sciences, Toronto, ON; Springer, Cham
(2018), 258 pp.

\medskip

\bibitem{EV21}
M. ~Efendiev, V. ~Vougalter. {\em Solvability in the sense of sequences
for some fourth order non-Fredholm operators,} J. Differential Equations,
{\bf 271} (2021), 280--300

\medskip

\bibitem{EV211}
M. ~Efendiev, V. ~Vougalter. {\em Existence of solutions for some non-Fredholm
integro-differential equations with mixed diffusion,} J. Differential
Equations, {\bf 284} (2021), 83--101

\medskip

\bibitem{EV25}
M. ~Efendiev, V. ~Vougalter. {\em On the solvability of some systems of
integro-differential equations with and without a drift,} Complex Var. Elliptic Equ.,
{\bf 70} (2025), no. 7, 1189--1207

\medskip

\bibitem{GS05}
H.G. ~Gebran, C.A. ~Stuart. {\em Fredholm and properness properties of
quasilinear elliptic systems of second order,} Proc. Edinb. Math. Soc. (2),
{\bf 48} (2005), no. 1, 91--124

\medskip

\bibitem{GS10}
H.G. ~Gebran, C.A. ~Stuart. {\em Exponential decay and Fredholm properties
in second-order quasilinear elliptic systems,} J. Differential Equations,
{\bf 249} (2010), no. 1, 94--117

\medskip

\bibitem{GVA06}
S. ~Genieys, V. ~Volpert, P. ~Auger.
{\em Pattern and waves for a model in population dynamics with nonlocal
consumption of resources,}  Math. Model. Nat. Phenom, {\bf 1} (2006), no. 1,
65--82

\medskip

\bibitem{K64} M.A. ~Krasnosel'skii. {\em Topological methods in the theory of
nonlinear integral equations.} The Macmillan Company, New York (1964), XI, 395 pp.

\medskip

\bibitem{LM}
J.-L. ~Lions, E. ~Magenes. { \em Probl\`emes aux limites non
homog\`enes et applications,} vol. 1. Dunod, Paris (1968)

\medskip

\bibitem{M22} A.B. ~Muravnik. {\em Elliptic equations with translations of general form in a half-space,}
Math. Notes, {\bf 111} (2022), no. 4, 587--594

\medskip

\bibitem{MZ23} A.B. ~Muravnik, N.V. ~Zaitseva. {\em Classical solutions of hyperbolic differential-difference equations with differently directed translations,} Lobachevskii J. Math., {\bf 44} (2023), no. 3, 920--925

\medskip

\bibitem{PY02} D.E. ~Pelinovsky, J. ~Yang. {\em A normal form for
nonlinear resonance of embedded solitons,}  Proc. R. Soc. Lond., Ser. A, Math.
Phys. Eng. Sci., {\bf 458} (2002), no. 2022, 1469--1497

\medskip

\bibitem{S92} A.L. ~Skubachevskii. {\em Boundary value problems for differential-difference equations with incommensurable shifts,}
Russ. Acad. Sci., Dokl., Math., {\bf 45} (1992), no. 3, 695--699

\medskip


\bibitem{TV}
S. ~Trofimchuk, V. ~Volpert.
{\em Traveling waves in delayed reaction-diffusion equations in biology,}
Math. Biosci. Eng., {\bf 17} (2020) no. 6, 6487--6514



\bibitem{Volevich}
L.R. ~Volevich. {\em Solubility of boundary value problems for
general elliptic systems,} Mat. Sb., {\bf 68} (1965),  no. 110, 373--416;
English translation: Am. Math. Soc. Transl., {\bf 67} (1968), no. 2,
182--225

\medskip

 \bibitem{V11}
 V.  ~Volpert. {\em Elliptic partial differential equations.}
 Volume 1:  Fredholm theory of elliptic problems in unbounded domains. Monographs in Mathematics, {\bf 101}.
 Basel: Birkh\"auser (2011), 639 pp.


\medskip

\bibitem{VKMP02} V. ~Volpert, B. ~Kazmierczak, M. ~Massot,
Z. ~Peradzynski. {\em Solvability conditions for elliptic problems with
non-Fredholm operators,} Appl. Math., {\bf 29} (2002), no. 2, 219--238

\medskip

\bibitem{VV1031} V. ~Volpert, V. ~Vougalter. {\em Solvability in the sense of
sequences to some non-Fredholm operators,} Electron. J. Differential Equations,
{\bf 2013}, Paper No. 160, 16 pp. (2013)

\medskip

\bibitem{VV08} V. ~Vougalter, V. ~Volpert. {\em Solvability conditions
for some non-Fredholm operators,} Proc. Edinb. Math. Soc. (2), {\bf 54}
(2011), no. 1,  249--271

\medskip

\bibitem{VV10} V. ~Vougalter, V. ~Volpert. {\em On the solvability
conditions for the diffusion equation with convection terms,}
Commun. Pure Appl. Anal., {\bf 11} (2012), no. 1, 365--373

\medskip

\bibitem{VV1021} V. ~Vougalter, V. ~Volpert. {\em On the existence
of stationary solutions for some non-Fredholm integro-differential
equations,} Doc. Math., {\bf 16} (2011), 561--580

\medskip

\bibitem{VV103} V. ~Vougalter, V. ~Volpert. {\em Solvability
conditions for some linear and nonlinear non-Fredholm elliptic
problems,} Anal. Math. Phys., {\bf 2} (2012), no. 4, 473--496

\medskip

\bibitem{VV25} V. ~Vougalter,  V. ~Volpert. {\em Solvability of some
integro-differential equations with the logarithmic Laplacian,} Appl. Anal.,
{\bf 104} (2025), no. 6, 1021--1035

\medskip


\bibitem{Z25} N.V. ~Zaitseva. {\em Classical solutions of hyperbolic equations with translation operators in lower-order derivatives,}
Differ. Equ., {\bf 61} (2025), no. 5, 627--649


\medskip

\bibitem{ZM23} N.V. ~Zaitseva, A.B. ~Muravnik. {\em A classical solution to a hyperbolic differential-difference equation with
a translation by an arbitrary vector,} Russ. Math., {\bf 67} (2023), no. 5, 29--34

\end{thebibliography}
\end{document}